\documentclass[12pt]{etds}            \usepackage{amsfonts,amssymb,latexsym,graphicx,psfrag} \usepackage{hyperref}        \newcommand{\D}{\displaystyle}
\newcommand{\restr}[1]{\raisebox{-0.3em}{$\lb|_{#1}\rb.$}} 
\newcommand{\ignore}[1]{}    
\newcommand{\breath}{\medskip} 
\newtheorem{thm}{Theorem}[section]

\newcounter{claimcount}[thm]  
\newcounter{XMPL}[thm]

\newtheorem{prop}[thm]{Proposition} 

\newtheorem{lemma}[thm]{Lemma} 

\newtheorem{cor}[thm]{Corollary}

\newcommand{\dfn}{\sf\em} 
\newcommand{\Theorem}[2]{\begin{thm}{\sf #1}  #2 \end{thm}}
\newcommand{\Proposition}[2]{\begin{prop}{\sf #1}  #2 \end{prop}}
\newcommand{\Lemma}[2]{\begin{lemma}{\sf #1}  #2 \end{lemma}}
\newcommand{\Corollary}[2]{\begin{cor}{\sf #1}  #2 \end{cor}} 
\newcommand{\thmfont}[1]{{\sl #1}}    
\newcommand{\examples}[1]{          \refstepcounter{thm}  \paragraph*{\sc Example \thethm:}  \begin{list}{$\langle$\alph{XMPL}$\rangle$}{\usecounter{XMPL}} 			{\setlength{\leftmargin}{0em} 			\setlength{\rightmargin}{0em}}   #1 	\hfill$\diamondsuit$\end{list}   			} 
\newcommand{\bthmlist}{ \begin{list}{{\bf(\alph{enumi})}} {\usecounter{enumi} \setlength{\leftmargin}{1em} \setlength{\itemsep}{0.2em} \setlength{\topsep}{0.2em} \setlength{\itemindent}{0em} \setlength{\parsep}{0em} \setlength{\rightmargin}{0em}} } 
\newcommand{\ethmlist}{\end{list}}    
\newcommand{\Claim}[1]{\refstepcounter{claimcount}                \noindent {\sc Claim \theclaimcount: \ }\thmfont{ #1}} 
\newcommand{\bprf}[1][Proof.]{\begin{list}{} 			{\setlength{\leftmargin}{0.7em} 			\setlength{\rightmargin}{0em} 			\setlength{\listparindent}{1em}}                         \item {\em \hspace{-1em}  #1  }} 
\newcommand{\eprf}{\end{list}} 
\newcommand{\bprfof}[1]{\begin{list}{} 			{\setlength{\leftmargin}{0.7em} 			\setlength{\rightmargin}{0em}}                         \item {\em \hspace{-1em}  Proof of #1.  }} 
\newcommand{\bthmprf}{\bprf}
\newcommand{\bclaimprf}{\bprf}
\newcommand{\ethmprf}{ \hfill$\Box$  \eprf  \breath  } 
\newcommand{\QED}{\hfill\ensuremath{\Box}}
\newcommand{\qed}{\QED}     
\newcommand{\beq}{\begin{eqnarray*}}
\newcommand{\eeq}{\end{eqnarray*}} 
\newcommand{\beqn}{ \begin{equation} }
\newcommand{\eeqn}{ \end{equation} }
\newcommand{\blist}{\begin{enumerate}}
\newcommand{\elist}{\end{enumerate}} 
\newcommand{\bdesc}{\begin{description}}
\newcommand{\edesc}{\end{description}}   
\newcommand{\Cesaro}{Ces\`aro }
\newcommand{\done}{{\mathsf{ 1\!\!1}}} 
\newcommand{\dB}{{\mathbb{B}}}
\newcommand{\dC}{{\mathbb{C}}}
\newcommand{\dE}{{\mathbb{E}}}
\newcommand{\dF}{{\mathbb{F}}}
\newcommand{\dG}{{\mathbb{G}}}
\newcommand{\dH}{{\mathbb{H}}}
\newcommand{\dI}{{\mathbb{I}}}
\newcommand{\dJ}{{\mathbb{J}}}
\newcommand{\dK}{{\mathbb{K}}}
\newcommand{\dL}{{\mathbb{L}}}
\newcommand{\dM}{{\mathbb{M}}}
\newcommand{\dN}{{\mathbb{N}}}
\newcommand{\dR}{{\mathbb{R}}}
\newcommand{\dS}{{\mathbb{S}}}
\newcommand{\dT}{{\mathbb{T}}}
\newcommand{\dU}{{\mathbb{U}}}
\newcommand{\dV}{{\mathbb{V}}}
\newcommand{\dW}{{\mathbb{W}}}
\newcommand{\dX}{{\mathbb{X}}}
\newcommand{\dZ}{{\mathbb{Z}}}       
\newcommand{\bA}{{\mathbf{ A}}}
\newcommand{\bB}{{\mathbf{ B}}}
\newcommand{\bP}{{\mathbf{ P}}}
\newcommand{\bU}{{\mathbf{ U}}}
\newcommand{\bV}{{\mathbf{ V}}}
\newcommand{\bX}{{\mathbf{ X}}}
\newcommand{\bY}{{\mathbf{ Y}}}
\newcommand{\bZ}{{\mathbf{ Z}}} 
\newcommand{\ba}{{\mathbf{ a}}}
\newcommand{\bb}{{\mathbf{ b}}}
\newcommand{\bc}{{\mathbf{ c}}}
\newcommand{\br}{{\mathbf{ r}}}
\newcommand{\bs}{{\mathbf{ s}}}
\newcommand{\bu}{{\mathbf{ u}}}
\newcommand{\bv}{{\mathbf{ v}}}
\newcommand{\bw}{{\mathbf{ w}}}
\newcommand{\bx}{{\mathbf{ x}}}
\newcommand{\brho }{{\boldsymbol{\rho}}}
\newcommand{\bxi}{{\boldsymbol{\xi }}}
\newcommand{\bchi}{{\boldsymbol{\chi }}}
\newcommand{\sA}{{\mathcal{ A}}}
\newcommand{\sB}{{\mathcal{ B}}}
\newcommand{\sF}{{\mathcal{ F}}}
\newcommand{\sI}{{\mathcal{ I}}}
\newcommand{\sM}{{\mathcal{ M}}}
\newcommand{\sP}{{\mathcal{ P}}}
\newcommand{\sX}{{\mathcal{ X}}}
\newcommand{\gA}{{\mathfrak{ A}}}
\newcommand{\gB}{{\mathfrak{ B}}}
\newcommand{\gC}{{\mathfrak{ C}}}
\newcommand{\gE}{{\mathfrak{ E}}}
\newcommand{\gF}{{\mathfrak{ F}}}
\newcommand{\gL}{{\mathfrak{ L}}}
\newcommand{\gO}{{\mathfrak{ O}}}
\newcommand{\gR}{{\mathfrak{ R}}}
\newcommand{\gS}{{\mathfrak{ S}}}
\newcommand{\gU}{{\mathfrak{ U}}}
\newcommand{\gW}{{\mathfrak{ W}}}
\newcommand{\gX}{{\mathfrak{ X}}}
\newcommand{\alp }{\alpha}
\newcommand{\bet }{\beta}
\newcommand{\gam }{\gamma}
\newcommand{\del }{\delta}
\newcommand{\eps }{\epsilon}
\newcommand{\sig }{\sigma} 
\newcommand{\omg }{\omega}
\newcommand{\Gam }{\Gamma}
\newcommand{\Omg }{\Omega}   
\newcommand{\h}{\widehat} 
\newcommand{\hsA}{{\widehat{\mathcal{ A}}}}
\newcommand{\fb}{{\mathsf{ b}}}
\newcommand{\ff}{{\mathsf{ f}}}
\newcommand{\fg}{{\mathsf{ g}}}
\newcommand{\fh}{{\mathsf h}}
\newcommand{\sfi}{{\mathsf{ i}}}
\newcommand{\fj}{{\mathsf{ j}}}
\newcommand{\fk}{{\mathsf{ k}}}
\newcommand{\fl}{{\mathsf{ l}}}
\newcommand{\fm}{{\mathsf{ m}}}
\newcommand{\fn}{{\mathsf{ n}}}
\newcommand{\fu}{{\mathsf{ u}}}
\newcommand{\fv}{{\mathsf{ v}}}
\newcommand{\fx}{{\mathsf{ x}}}
\newcommand{\tl}{\widetilde} 
\newcommand{\tlM}{{\widetilde{M}}}
\newcommand{\tlm}{{\widetilde{m}}}
\newcommand{\tlw}{{\widetilde{w}}}
\newcommand{\lb}{\left}
\newcommand{\rb}{\right} 
\newcommand{\maketall}{\rule[-0.5em]{0em}{1em}}        
\newcommand{\map}{{\longrightarrow}}
\newcommand{\goto}{{\rightarrow}}
\newcommand{\into}{{\map}}
\newcommand{\image}[1]{\mathrm{image}\lb(#1\rb)}  
\newcommand{\seilpmi}{{\Longleftarrow}}
\newcommand{\statement}[1]{\lb(  \maketall       \begin{minipage}{40em}       \begin{tabbing}         #1        \end{tabbing}      \end{minipage}  \rb)}     
\newcommand{\oo}{{\infty}}        
\newcommand{\tensor}{\otimes}
\newcommand{\Tensor}{\bigotimes}
\newcommand{\union}{\cup}
\newcommand{\Union}{\bigcup}
\newcommand{\intsct}{\cap}
\newcommand{\Intsct}{\bigcap}
\newcommand{\disj}{\sqcup}
\newcommand{\Disj}{\bigsqcup}   
\newcommand{\set}[2]{{\left\{ #1 \; ; \; #2 \right\} }} 
\newcommand{\supp}[1]{{\sf supp}\lb(#1\rb)}      
\newcommand{\inn}[1]{{\left\langle #1 \right\rangle }}       
\newcommand{\chr}[1]{{{\done}_{{#1}}}} 
\newcommand{\choice}[1]{{\lb\{ \begin{array}{rcl}                                 #1                                \end{array}  \rb.  }}                     %
\newcommand{\eeequals}[1]{\raisebox{-0.9ex}{$\overline{\overline{{\scriptscriptstyle{\mathrm{#1}}}}}$}} 
\newcommand{\leeeq}[1]{\raisebox{-1ex}{${{\D\leq} \atop {\scriptscriptstyle{\mathrm{#1}}}}$}} 
\newcommand{\lt}[1]{\raisebox{-1ex}{${{\D<} \atop {\scriptscriptstyle{\mathrm{#1}}}}$}} 
\newcommand{\grt}[1]{\raisebox{-1ex}{${{\D>} \atop {\scriptscriptstyle{\mathrm{#1}}}}$}}  
\newcommand{\geeeq}[1]{\raisebox{-1ex}{${{\geq} \atop {\scriptscriptstyle{\mathrm{#1}}}}$}} 
\newcommand{\closeto}[1]{{{\raisebox{-1ex}    {$\widetilde{\ {\scriptstyle #1}\ }$}}}}      
\newcommand{\cl}[1]{{\sf cl}\lb(#1\rb)}
\newcommand{\shift}[1]{\sig^{#1}}    
\newcommand{\Expct}[2][]{\dE_{#1}\lb[#2\rb]}
\newcommand{\Prob}[1]{\ensuremath{{\sf Prob}\lb[ #1 \rb] }}
\newcommand{\goesto}[2]{{ -\!\!\!-\!\!\!-\!\!\!-\!\!\!\!\!\!\!\!\!\!\!  ^{{\scriptscriptstyle #2}}_{{\scriptscriptstyle #1}}   \!\!\!\!\!\!\!\!\!\longrightarrow }}                         
\newcommand{\wkstlim}{\mathrm{wk}^*\!\!-\!\!\!\lim}       
\newcommand{\card}[1]{\#\lb(#1\rb)}
\newcommand{\diam}[1]{{\sf diam}\lb[#1\rb]} 
\newcommand{\Meas}[1]{{\sM\lb[#1\rb] }}
\newcommand{\mtrx}[3]{{\lb[#1  |_{#2}^{#3} \rb]}}         
\newcommand{\Matrix}[1]{{\lb[\begin{array}{ccccccccccccccccccccccccr} #1 \end{array}\rb]}}   
\newcommand{\rank}[1]{{{\sf rank}\lb[#1\rb]}}
\newcommand{\Real}{\dR}
\newcommand{\Natur}{\dN}
\newcommand{\Zahl}{\dZ}
\newcommand{\Zahlmod}[1]{{\Zahl_{/#1}}}
\newcommand{\Cplx}{\dC}
\newcommand{\Torus}[1]{{{\dT}^{#1}}} 
\newcommand{\CC}[1]{{\lb[ #1 \rb]}}
\newcommand{\CO}[1]{{\lb[ #1 \rb)}}
\newcommand{\OC}[1]{{\lb( #1 \rb]}}
\newcommand{\OO}[1]{{\lb( #1 \rb)}}   
\renewcommand{\implies}{\ensuremath{\Longrightarrow}}\renewcommand{\And}{\mbox{\ and \ }} 
\newcommand{\Prec}[1]{\sP\lb(#1\rb)}
\newcommand{\Fol}[1]{\sF\lb(#1\rb)}
\newcommand{\bound}{\partial\,} 
\newcommand{\cyl}[1]{\lb[#1\rb]}  
\newcommand{\Selfmap}[1]{#1\into#1}   
\newcommand{\Lucas}[1]{\dL\lb(#1\rb)}   
\newcommand{\density}[1]{{\sf density}\lb(#1\rb)}   
\newcommand{\reldense}[1]{{\sf rel \, density}\lb[#1\rb]}  
\newcommand{\comesfrom}{\dashleftarrow}
\newcommand{\goesinto}{\dashrightarrow} 
\newcommand{\Lat}{\dM}
\newcommand{\Ball}{\dB} 
\newcommand{\Srank}[2][S]{{\sf rank}_{#1}\lb(#2\rb)}          
\newcommand{\centre}[1]{{\sf centre}\lb(#1\rb)} 

\newcommand{\SFT}{\gF}
\newcommand{\Sof}{\gS}
\newcommand{\bsof}{\bs}
\newcommand{\bsft}{\bx}
\newcommand{\sof}{s}
\newcommand{\lat}{\fm}
\newcommand{\ecp}{\hfill~$\Diamond$} 
\newcommand{\Even}{\gE}
\newcommand{\Odd}{\gO}   
\newcommand{\dm}[1]{\lb|\!\lb[#1\rb]\!\rb|}
\newcommand{\ldm}[1]{{\lb\langle\!\lb\langle#1\rb\rangle\!\rb\rangle}} 
\newcommand{\dilate}[2]{#1\circ\Phi^{#2\cdot \ldm{#1}}} 
\newcommand{\logdiam}[1]{\lb\lceil\maketall \log_p\dm{#1}\rb\rceil }  
\newcommand{\tldJ}{\widetilde{\dJ}}
\newcommand{\tldH}{\widetilde{\dH}} 
\newcommand{\eclaimprf}{ \hfill $\Diamond$~{\scriptsize {\tt Claim~\theclaimcount}}\eprf} 
\newcommand{\THET}{\Theta} 

\begin{document}

\ETDS{1}{26}{26}{2006} 
\title{Asymptotic Randomization of Sofic Shifts by Linear Cellular Automata\footnote{This research was partially supported by NSERC Canada, and 
was also supported
by the kind hospitality of  the Universidad de Chile during July 2003.}}

\runningheads{M. Pivato and R. Yassawi}{Randomization of Sofic Shifts by Linear Cellular Automata}

\author{Marcus Pivato \ and \ Reem Yassawi}

\address{   Department of Mathematics, Trent University, 1600 West Bank Drive,
Peterborough, Ontario,  K9J 7B8,  Canada
 \\ \email{{\tt marcuspivato@trentu.ca} and {\tt ryassawi@trentu.ca}} }

\recd{July 31, 2004;  \ accepted April 6, 2006}


\begin{abstract}
  Let $\Lat=\Zahl^D$ be a $D$-dimensional lattice, and let $(\sA,+)$
be an abelian group.  $\sA^\Lat$ is then a compact abelian group under
componentwise addition.  A continuous function $\Phi:\sA^\Lat\into\sA^\Lat$ 
is called a {\em linear cellular automaton} (LCA) if
there is a finite subset $\dF\subset\Lat$ and nonzero coefficients
$\varphi_\ff\in\Zahl$ so that, for any $\ba\in\sA^\Lat$, \quad $
\Phi(\ba) \ = \ \sum_{\ff\in\dF} \varphi_\ff\cdot\shift{\ff}(\ba)$.

 Suppose $\mu$ is a probability measure on $\sA^\Lat$ whose support is
a subshift of finite type or sofic shift.  We provide sufficient
conditions (on $\Phi$ and $\mu$) under which $\Phi$ {\em
asymptotically randomizes} $\mu$, meaning that $\D\wkstlim_{\dJ\ni
j\goto\oo} \Phi^j\mu \ = \ \eta$, where $\eta$ is the Haar measure on
$\sA^\Lat$, and $\dJ\subset\Natur$ has \Cesaro density 1.  In the case
when $\Phi=1+\shift{}$ and $\sA=(\Zahlmod{p})^s$ ($p$ prime), we
provide a condition on $\mu$ that is both necessary and sufficient.
We then use this to construct zero-entropy measures which are randomized
by $1+\shift{}$.

\end{abstract}

  Let $D\geq 1$, and let $\Lat:=\Zahl^D$ be the $D$-dimensional lattice.
If $\sA$ is a (discretely topologised) finite set, then $\sA^\Lat$ is
compact in the Tychonoff topology. 
For any $\fv\in\Lat$, let $\shift{\fv}:\Selfmap{\sA^\Lat}$ be the shift
map: $\shift{\fv}(\ba) \ := \ \mtrx{b_\lat}{\lat\in\Lat}{}$, where
$b_\lat := a_{\lat-\fv}$, \ $\forall\lat\in\Lat$.  
 A {\dfn cellular automaton} (CA) is a continuous map
$\Phi:\Selfmap{\sA^\Lat}$ which commutes with all shifts: \ for any
$\lat\in\Lat$,\ \ \ $\shift{\lat}\circ
\Phi \ = \ \Phi \circ \shift{\lat}$.
Let $\eta$ be the uniform Bernoulli measure on $\sA^\Lat$.  If $\mu$ is
 another probability measure on $\sA^\Lat$, we say  $\Phi$ {\dfn
 asymptotically randomizes} $\mu$ if $\D\wkstlim_{\dJ\ni j \goto\oo}
\Phi^j \mu \ = \ \eta$, where $\dJ\subset\Natur$ has \Cesaro
density one.

  If $(\sA,+)$ is a finite abelian group, then $\sA^\Lat$ is a product
group, and $\eta$ is the Haar measure. A {\dfn linear
cellular automaton} (LCA) is a CA $\Phi$ with a finite subset
$\dF\subset\Lat$ (with $\card{\dF}\geq 2$),
and nonzero coefficients $\varphi_\ff\in\Zahl$ (for all $\ff\in\dF$) so
that, for any $\ba\in\sA^\Lat$, 
\beqn
\label{LCA}
\Phi(\ba) \quad = \quad \sum_{\ff\in\dF} \varphi_\ff\cdot\shift{\ff}(\ba).
\eeqn
Linear cellular automata are known to asymptotically randomize a wide
variety of measures
\cite{MaassMartinez,MaassMartinezII,MaassHostMartinez,Lind,FerMaassMartNey},
including those satisfying a correlation-decay condition called {\em harmonic mixing}
\cite{PivatoYassawi1,PivatoYassawi2,MaassMartinezPivatoYassawi}.  However, all known 
sufficient conditions for asymptotic randomization  (and
for harmonic mixing, in particular) require $\mu$ to have {\dfn full
support}, i.e.  $\supp{\mu}= \sA^\Lat$.
 
 We here investigate asymptotic randomization when $\supp{\mu}\subsetneq
 \sA^\Lat$.  In particular we consider the case when $\supp{\mu}$ is a
 sofic shift or subshift of finite type.  In
\S\ref{S:mrf}, we let $\sA=\Zahlmod{p}$ ($p$ prime), and
demonstrate asymptotic randomization for any Markov random field that
is {\em locally free}, a much weaker assumption than full support.
However, in \S\ref{S:non.harm.mixing} we show that harmonic mixing is
a rather restrictive condition, by exhibiting a measure whose support
is a mixing sofic shift but which is {\em not} harmonically mixing.

Thus, in \S\ref{S:wkhm}, we introduce the less restrictive concept of
{\em dispersion mixing} (for measures) and the dual concept of {\em
dispersion} (for LCA), and state our main result: any dispersive LCA
asymptotically randomizes any dispersion mixing measure.  In
\S\ref{S:disperse}, we let $\sA=(\Zahlmod{p})^s$ ($p$ prime, $s\in\Natur$)
and introduce {\em bipartite} LCA, a broad class
exemplified by the automaton $1+\shift{}$.  We then show that any
bipartite LCA is dispersive.

In \S\ref{S:unimix}, we show that any {\em uniformly mixing} and {\em
harmonically bounded} measure is dispersion mixing.  In particular, in
\S\ref{S:mkwrd}, we show this implies that any mixing Markov measure
(supported on a subshift of finite type), and any continuous factor of
a mixing Markov measure (supported on a sofic shift) is dispersion
mixing, and thus, is asymptotically randomized by any dispersive LCA
(e.g. $1+\shift{}$).  Thus, the example of
\S\ref{S:non.harm.mixing} {\em is} asymptotically randomized, even
though it is not harmonically mixing.  

  In \S\ref{S:lucas.mix}, we refine the results of 
\S\ref{S:wkhm}-\S\ref{S:disperse} by introducing {\em Lucas mixing},
(a weaker condition than dispersion mixing).  When
$\sA=(\Zahlmod{p})^s$, we show that a
measure is asymptotically randomized by the automaton
$1+\shift{}$ if and {\em only if} it is Lucas mixing.
Finally, in \S\ref{S:zero.random}, we use Lucas mixing to construct a class of
{\em zero-entropy} measures which are asymptotically randomized by
randomized by $1+\shift{}$, thereby refuting the conjecture that positive
entropy is necessary for asymptotic randomization. 

\subsubsection*{Preliminaries \& Notation:}

  Throughout, $(\sA,+)$ is an abelian group (usually
$\sA=(\Zahlmod{p})^s$,where $p$ is prime and $s\in\Natur$).  Elements
of $\sA^\Lat$ are denoted by boldfaced letters (e.g.  $\ba$, $\bb$,
$\bc$), and subsets by gothic letters (e.g. $\gA$, $\gB$, $\gC$).
Elements of $\Lat$ are sans serif (e.g. $\fl$, $\fm$, $\fn$) and
subsets are $\dU,\dV,\dW$.

If  $\dU\subset\Lat$ and $\ba\in\sA^\Lat$ then $\ba_\dU :=
\mtrx{a_\fu}{\fu\in\dU}{}$ is the `restriction' of $\ba$ to an element of
$\sA^\dU$.  For any  $\bb\in\sA^\dU$, let
$\cyl{\bb} \ := \ \set{\bc\in\sA^\Lat}{\bc_\dU=\bb}$ be the corresponding
cylinder set.  In particular, if $\ba\in\sA^\Lat$, then
$\cyl{\ba_\dU} \ := \ \set{\bc\in\sA^\Lat}{\bc_\dU=\ba_\dU}$.

\paragraph*{Measures:}
  Let $\sM(\sA^\Lat)$ be the set of Borel probability measures on
$\sA^\Lat$.  If $\mu\in\sM(\sA^\Lat)$ and $\dI\subset\Lat$, then
let $\mu_\dI\in\sM(\sA^\dI)$ be the marginal projection of $\mu$
onto $\sA^\dI$.  If $\dJ\subset\Lat$ and $\bb\in\sA^\dJ$, then
let $\mu^{(\bb)}\in\sM(\sA^\Lat)$ be the conditional probability
measure in the cylinder set $\cyl{\bb}$.  
In other words, for any $\gX\subset\sA^\Lat$,
\quad $\mu^{(\bb)}[\gX] \ := \ 
\mu\lb(\gX\intsct\cyl{\bb}\rb)/\mu\cyl{\bb}$.
  In particular, if $\dI\subset\Lat$ is finite, then
$\mu^{(\bb)}_\dI\in\sM(\sA^\Lat)$ is the conditional probability
measure on the $\dI$ coordinates: \  for any $\bc\in\sA^\dI$,
\quad $\D\mu_\dI^{(\bb)}[\bc] \ := \
\mu\lb(\cyl{\bc}\intsct\cyl{\bb}\rb)/\mu\cyl{\bb}$.

\paragraph{Subshifts:}
    A {\dfn subshift} \cite{Kitchens,LindMarcus} is a a closed,
 shift-invariant subset $\gX\subset\sA^\Lat$.  If $\dU\subset\Lat$,
 then let $\gX_\dU:=\set{\bx_\dU}{\bx\in\gX}$ be all {\dfn admissible
 $\dU$-blocks} in $\gX$.  If $\dU\subset\Lat$ is finite, and $\gW =
 \{\bw_1,\ldots,\bw_N\}\subset
\sA^\dU$ is a collection of admissible blocks, then the
induced {\dfn subshift of finite type} (SFT) is the largest
subshift $\gX\subset\sA^\Lat$ such that $\gX_\dU=\gW$. 
In other
words, $ \gX \ := \ \Intsct_{\lat\in\Lat} \shift{\lat}\cyl{\gW}$,
where $\cyl{\gW} \ := \  \set{\ba\in\sA^\Lat}{\ba_\dU\in\gW}$.
A {\dfn sofic shift} is the image of an SFT under a block map.

  In particular, if $\Lat=\Zahl$ and $\dU=\{0,1\}$, then $\gX$
is called  {\dfn topological Markov shift}, and the {\dfn transition
matrix} of $\gX$ is the matrix $\bP=[p_{ab}]_{a,b\in\sA}$, where
$p_{ab}=1$ if $[ab]\in\gW$, and $p_{ab}=0$ if $[ab]\not\in\gW$.

\paragraph*{Characters:}
Let $\Torus{1}\subset\Cplx$ be the circle group.
 A {\dfn character} of $\sA^\Lat$ is a continuous homomorphism
$\bchi:\sA^\Lat\into\Torus{1}$; the group of such characters is denoted
$\h{\sA^\Lat}$.  For any $\bchi\in\h{\sA^\Lat}$ there is a finite
subset $\dK\subset\Lat$, and nontrivial $\chi_\fk\in\hsA$ for all
$\fk\in\dK$, such that, for any $\ba\in\sA^\Lat$, \quad $\bchi(\ba) \ = \
\D\prod_{\fk\in\dK} \chi_\fk(a_\fk)$.  We indicate this by writing:
``$\D \bchi\ =\ \Tensor_{\fk\in\dK} \chi_\fk$''.
The {\dfn rank} of $\bchi$ is the cardinality of $\dK$.

\paragraph{\Cesaro Density:}

  If $\ell,n\in\Zahl$, then let $\CO{\ell...n}:=\set{m\in\Zahl}{\ell\leq m < n}$.
If $\dJ\subset\Natur$, then the {\dfn \Cesaro density} of $\dJ$ is
defined: \ $\density{\dJ} \ := \ \D\lim_{N\goto\oo} \frac{1}{N}
\card{\maketall \dJ\intsct\CO{0...N}}$.  If $\dJ,\dK\subset\Natur$, then their
{\dfn relative \Cesaro density} is defined: \
\[
\reldense{\dJ/\dK} \quad := \quad
\lim_{N\goto\oo} \frac{\card{\dJ\intsct\CO{0...N}}}
			{\card{\dK\intsct\CO{0...N}}}.
\]
In particular, $\density{\dJ} \ = \ \reldense{\dJ/\Natur}$.

\section{\label{S:mrf}
Harmonic Mixing of Markov Random Fields}

   Let $\Ball\subset\Lat$ be a finite subset, symmetric under
multiplication by $-1$ (usually, 
$\Ball=\{-1,0,1\}^D$).  For any $\dU\subset\Lat$, we define
\[
  \cl{\dU} \ := \ \set{\fu+\fb}{\fu\in\dU \And \fb\in\Ball}
\quad\quad\And\quad\quad
\partial\dU \ := \ \cl{\dU} \setminus \dU.
\]
For example, if $\Lat=\Zahl$ and $\Ball=\{-1,0,1\}$, then
$\partial\{0\} \ = \ \{\pm1\}$.

  Let $\mu\in\sM(\sA^\Lat)$.  Suppose $\dU\subset\Lat$, and let
$\dV:=\partial\dU$ and $\dW = \Lat\setminus\cl{\dU}$.  If
$\bb\in\sA^\dV$, then we say $\bb$ {\dfn isolates} $\dU$ from $\dW$ if
the conditional measure $\mu^{(\bb)}$ is a product of
$\mu^{(\bb)}_\dU$ and $\mu^{(\bb)}_\dW$.  That is, for any
$\gU\subset\sA^\dU$ and $\gW\subset\sA^\dW$, we have
$\mu^{(\bb)}\lb(\gU\intsct\gW\rb) \ = \
\mu^{(\bb)}_\dU(\gU) \cdot \mu^{(\bb)}_\dW(\gW)$.

  We say that $\mu$ is a {\dfn Markov random field}
 \cite{Bremaud,KindermannSnell} with {\dfn interaction range} $\Ball$
 (or write, ``$\mu$ is a {\dfn $\Ball$-MRF}'') if, for any
 $\dU\subset\Lat$ with $\dV=\partial\dU$ and $\dW =
 \Lat\setminus\cl{\dU}$, any choice of $\bb\in\sA^\dV$ isolates $\dU$
 from $\dW$.

For example, if $\Lat=\Zahl$ and $\Ball=\{-1,0,1\}$, then $\mu$ is a
$\Ball$-MRF iff $\mu$ is a (one-step) Markov chain.  If
$\Ball=\CC{-N...N}$, then $\mu$ is a $\Ball$-MRF iff $\mu$ is an $N$-step
Markov chain.

\Lemma{}
{ If $\mu$ is a Markov random field, then $\supp{\mu}$ is a subshift of finite type.\qed}

  For example, if $\mu$ is a Markov chain on $\sA^\Zahl$, then $\supp{\mu}$
is a topological Markov shift.

\breath

  Let $\Ball\subset\Lat$, and let $\mu\in\sM(\sA^\Lat)$ be $\Ball$-MRF.  Let
$\dS:=\Ball\setminus\{0\}$.  For any $\bb\in\sA^\dS$, let
$\mu_0^{(\bb)}\in\sM(\sA)$ be the conditional probability measure on
the zeroth coordinate.  We say that $\mu$ is {\dfn locally free} if,
for any $\bb\in\sA^\dS$, \quad $\card{\supp{\mu_0^{(\bb)}}}\geq 2$.

\paragraph*{\sc Example:} 
If $D=1$, then $\Ball=\{-1,0,1\}$,
\ $\dS=\{\pm1\}$, and $\mu$ is a Markov chain.  Thus, $\supp{\mu}$
is a topological Markov shift, with transition matrix
$\bP=[p_{ab}]_{a,b\in\sA}$.  For any $a,b\in\sA$, write $a\leadsto b$
if $p_{ab}=1$, and define the {\dfn follower} and {\dfn predecessor}
sets
\[
 \Fol{a}\ :=\quad\set{b\in\sA}{a\leadsto b}
\ \And \ 
 \Prec{b}\ :=\quad\set{a\in\sA}{a\leadsto b}.
\]
It is easy to show that the following are equivalent:

\blist
\item $\mu$ is locally free.

\item  Every entry of $\bP^2$ is 2 or larger.

\item For any $a,b\in\sA$, $\card{\Fol{a}\intsct\Prec{b}}\geq 2$.
\elist

 Recall that $\hsA$ is the dual group of $\sA$.
For any $\bchi\in\hsA$ and $\nu\in\sM(\sA)$, let
$\D \inn{\chi,\nu} \ := \  \sum_{a\in\sA} \chi(a)\cdot\nu\{a\}$.
It is easy to check:

\Lemma{\label{cond.exp.lem}}
{   Let $p$ be prime and $\sA=\Zahlmod{p}$.
If $\mu$ is a locally free MRF on $\sA^\Lat$, then there is some
  $c<1$ such that, for all nontrivial $\chi\in\hsA$, and any $\bb\in\sA^\dS$,
\ \ $\lb| \inn{\chi, \mu^{(\bb)}_0} \rb| \quad\leq\quad c$.
\qed
}
For any  $\bchi\in\h{\sA^\Lat}$ and $\mu\in\sM\lb(\sA^\Lat\rb)$, define
\ $\D
\inn{\bchi,\mu} \ := \ \int_{\sA^\Lat} \bchi(\ba)\ d\mu[\ba]$. \ 
  A measure $\mu$ is called {\dfn harmonically mixing} if,
for any $\eps>0$, there is some $R\in\Natur$ such that,
for any  $\bchi\in\h{\sA^\Lat}$,
\[
\statement{$\rank{\bchi} \ > \ R$}
\implies
\statement{$\lb|\inn{\bchi,\mu}\rb| \ < \ \eps$}.
\]
The significance of this is the following \cite[Theorem 12]{PivatoYassawi1}:

\paragraph*{\sc Theorem:}\ 
{\sl  Let $\sA=\Zahlmod{p}$,  where $p$ is prime.
Any  LCA on $\sA^\Lat$  asymptotically randomizes any
harmonically mixing measure.}\qed

\breath

  Most MRFs with full support are harmonically mixing \cite[Theorem 15]{PivatoYassawi2}.  We now extend this.

\Theorem{\label{mrf.local.free.HM}}
{
 Let $\sA=\Zahlmod{p}$, where $p$ is prime.
 Any locally free MRF on $\sA^\dM$ is harmonically mixing.
}
\bthmprf Let $\mu$ be a locally free $\dB$-MRF.
  A subset $\dI\subset\Lat$ is {\dfn $\dB$-separated} if $(\sfi-\fj) \not\in \Ball$
for all $\sfi,\fj\in\dI$ with $\sfi\neq \fj$.
  Let $\dK\subset\Lat$ be  finite, and let
 $\bchi\ := \ \D\Tensor_{\fk\in\dK}\chi_\fk$ be a character
of $\sA^\Lat$.

\Claim{Let $K:=\card{\dK}=\rank{\bchi}$, and let 
$B \ := \ \max \set{|\fb_1-\fb_2|}{\fb_1,\fb_2\in\Ball}$.
There exists a $\dB$-separated subset $\dI\subset\dK$ such that
\beqn
\label{I.bound}
\card{\dI} \quad=\quad I \quad \geq \quad \frac{K}{B^D}.
\eeqn
}
\bclaimprf  Let $\tl\Ball := \CO{0...B}^D$ be a box of sidelength $B$.  Cover
$\dK$ with disjoint translated copies of $\tl\Ball$, so that
\[
  \dK \quad \subset\quad \Disj_{\sfi\in\dI}\  \lb(\tl\Ball+\sfi\rb)
\]
 for some set $\dI\subset\dK$.  Thus, $|\sfi-\fj|\geq B$ for any $\sfi,\fj\in\dI$
 with $\sfi\neq \fj$, so $(\sfi-\fj) \not\in \Ball$.  Also, $\card{\tl\Ball} =
 B^D$, so each copy covers at most $B^D$ points in $\dK$.  Thus, we
 require at least $\frac{K}{B^D}$ copies to cover all of $\dK$.  In
 other words, $I\geq \frac{K}{B^D}$.
\eclaimprf

 Thus, \ $\bchi \ = \  \bchi_\dI\cdot \bchi_{\dK\setminus\dI}$,  \ where \
$\D
  \bchi_\dI(\ba) \ := \
\prod_{\sfi\in\dI} \chi_\sfi(a_\sfi)$, \ 
 and \
$\D  \bchi_{\dK\setminus\dI}(\ba) \ := \ 
\prod_{\fk\in\dK\setminus\dI} \chi_\fk(a_\fk)$.

Let $\dJ \ := \ (\bound\dI) \union (\dK\setminus\dI)$;
fix $\bb\in\sA^\dJ$, and let $\mu_\dI^{(\bb)}\in\sM(\sA^\dI)$ be the
corresponding conditional probability measure.  Since $\mu$ is a Markov
random field, and the $\dI$ coordinates are `isolated' from one another
by $\dJ$ coordinates, it follows that $\mu_\dI^{(\bb)}$ is a product measure.  In other words,
for any $\ba\in\sA^\dI$,
\beqn
\label{conditional.independence}
  \mu_\dI^{(\bb)}[\ba] \quad=\quad \prod_{\sfi\in\dI} \mu_\sfi^{(\bb)}\{a_\sfi\}.
\eeqn
Thus,
the conditional expectation of $\bchi_\dI$ is given:
\beq
\inn{\bchi_\dI, \ \mu_\dI^{(\bb)}}
&=&
\sum_{\ba\in\sA^\dI} \mu^{(\bb)}_\dI[\ba] \cdot 
\lb(\prod_{\sfi\in\dI} \chi_\sfi(a_\sfi)\rb)
\quad\eeequals{(*)} \ \quad
\sum_{\ba\in\sA^\dI} 
 \lb(\prod_{\sfi\in\dI} \mu^{(\bb)}_\sfi\{a_\sfi\}\cdot\chi_\sfi(a_\sfi)\rb)
\\ 
&=&
\prod_{\sfi\in\dI} \lb( \sum_{a_\sfi\in\sA} \mu^{(\bb)}\{a_\sfi\}\cdot \chi_\sfi(a_\sfi)\rb)
\quad=\quad
\prod_{\sfi\in\dI} \inn{\chi_\sfi, \ \mu^{(\bb)}_\sfi},
\eeq
where $(*)$ is by equation (\ref{conditional.independence}).
Thus, $\D
\inn{\bchi, \ \mu^{(\bb)}}
 = 
\bchi_{\dK\setminus\dI}(\bb)\cdot
\inn{\bchi_\dI, \ \mu_\dI^{(\bb)}}
 =
\bchi_{\dK\setminus\dI}(\bb)\cdot
\prod_{\sfi\in\dI} \inn{\chi_\sfi, \ \mu^{(\bb)}_\sfi}$.  
Thus, if $I=\card{\dI}$, then
\beqn
\label{cond.exp.bound}
  \lb|\inn{\bchi, \ \mu^{(\bb)}}\rb|
\quad=\quad
\lb|\maketall \bchi_{\dK\setminus\dI}(\bb)\rb|\cdot
\prod_{\sfi\in\dI} \lb|\inn{\chi_\sfi, \ \mu^{(\bb)}_\sfi}\rb|
\quad\leq\quad
1\cdot c^I
\eeqn
where the last step follows from Lemma \ref{cond.exp.lem}.
But \ $\D
\inn{\bchi,\mu}
\ = \ 
\sum_{\bb\in\sA^\dJ} \mu[\bb]\cdot\inn{\bchi, \ \mu^{(\bb)}}$, \
so
\[ \lb|\inn{\bchi,\mu}\rb|
\quad\leq\quad
\sum_{\bb\in\sA^\dJ} \mu[\bb]\cdot\lb|\inn{\bchi, \ \mu^{(\bb)}}\rb|
\quad \leeeq{(*)} \quad
\sum_{\bb\in\sA^\dJ} \mu[\bb]\cdot c^I
\quad = \quad  
c^I
\quad \leeeq{(\dagger)} \quad  c^{K/(B^D)} \goesto{K\goto\oo}{} 0.
\]
Here $(*)$ is by equation (\ref{cond.exp.bound}) and $(\dagger)$ is 
by equation (\ref{I.bound}).
\ethmprf

\section{The Even Shift is Not Harmonically Mixing
\label{S:non.harm.mixing}}

  We will now construct a measure $\nu$, supported on a sofic shift,
which is {\em not} harmonically mixing.  Nonetheless, we'll show in
\S\ref{S:wkhm}-\S\ref{S:unimix} that this measure {\em is} asymptotically randomized by many LCA.

Let $\gX \subset \lb(\Zahlmod{3}\rb)^\Zahl$ be the subshift of finite
type defined by the transition matrix
\[ \bA \quad=\quad \Matrix{
1 & 0 & 1 \\
1 & 0 & 1 \\
0 & 1 &0 },\ \ 
\mbox{where,} \ \forall i,j\in\Zahlmod{3}, \ 
\ a_{ij} = \choice{ 1 &&\mbox{if $j\leadsto i$ is allowed}\\
			0  &&\mbox{if $j\leadsto i$ is not allowed}}
\]
Let $\Phi: \gX \rightarrow \lb(\Zahlmod{2}\rb)^\Zahl$ be the factor map
of radius 0 
which sends $0$ into $0$ and both $1$ and $2$  to $1$.
 Then  $\Sof: =   \Phi(\gX)$ is Weiss's {\em Even Sofic Shift}: \  if 
$\bsof \, \in \Sof$, 
then there are an even number of 1's between any two occurrences of $0$ in
 $\bsof$.

 For any $N\in\Natur$, and $ i, j \in\Zahlmod{3}$, let 
$\gX_{ij}^{N} \ := \ \set{ \bsft \, \in \gX}{ x_{0} = i, \  x_{N} = j }$,
and let:
\[
 \Even_{N} \ :=\ 
  \set{ \bsof \, \in \Sof}{ \sum_{n=0}^{N} \sof_{n} \mbox{ is even }}, \ \mbox{and} \
  \Odd_{N}\ := \ 
 \set{ \bsof \, \in \Sof}{ \sum_{n=0}^{N} \sof_{n} \mbox{ is odd } }.
\] 

\Lemma{\label{even.odd}}
{$\forall i,j\in\Zahlmod{3}$, either
$\Phi \lb(\gX_{i,j}^{N}\rb) \subset \Even_{N}$ or $\Phi \lb(\gX_{i,j}^{N}\rb) \subset \Odd_{N}$.  In particular,
\ignore{\beq
 \Phi \lb( \gX_{0,0}^{N} \disj \gX_{1,2}^{N} \disj \gX_{2,1}^{N} \disj \gX_{0,2}^{N} \disj \gX_{1,0}^{N} \rb) &= &\Even_{N},\\
\mbox{and} \ \  
 \Phi \lb( \gX_{1,1}^{N} \disj \gX_{0,1}^{N} \disj \gX_{2,0}^{N} \disj \gX_{2,2}^{N}\rb)
&= &\Odd_{N}.
\eeq}
\[
\Phi \lb( \gX_{0,0}^{N} \disj \gX_{1,2}^{N} \disj \gX_{2,1}^{N} \disj \gX_{0,2}^{N} \disj \gX_{1,0}^{N} \rb) \ \ = \ \ \Even_{N}\qquad
\mbox{and} \qquad  
 \Phi \lb( \gX_{1,1}^{N} \disj \gX_{0,1}^{N} \disj \gX_{2,0}^{N} \disj \gX_{2,2}^{N}\rb)
\ \ = \ \ \Odd_{N}.
\]
}
\bthmprf  Let $\bsft \, \in
\gX_{ij}^{N}$, and $\bsof := \Phi(\bsft)$.   Note that, if $k<k^*$ are
any two values such that $x_k = 0 = x_{k^*}$, then $\D \sum_{n=k}^{k^*} \sof_n$
is even.  In particular, let $k$ be the first
element of $\CC{0...N}$ where $x_{k} =0$, and let $k^{*}$ be the
last element of $\CC{0...N}$  where  $x_{k^{*}}=0$.  
Thus, $\D\sum_{n= k}^{k^{*}} \sof_{n} \equiv 0 \pmod{2}$, so that
$\D\sum_{n=0}^N \sof_n \ \equiv \ \sum_{n=0}^{k-1} \sof_n + \sum_{n=k^*+1}^N \sof_n \pmod{2}$. 

  But since $x_{k-1} \neq 0 \neq x_{k^{*}+1}$ by construction, the
definition of $\gX$ forces $x_{k-1}= 2$ and $x_{k^{*}+1}=1$.
Thus the parity of $\D \sum_{n=0}^{k-1} \sof_{n}$ depends
only on the value of $x_{0} = i$.  Similarly the parity of
$\D\sum_{n=k^{*}+1}^{N} \sof_{n}$ depends only on $x_{N} = j$.
\ethmprf

  Let $\mu\in\Meas{\gX}$ be a mixing Markov measure on $\gX$,
 with transition matrix $\bP$ and Perron measure $\brho =
(\rho_0,\rho_1,\rho_2)\in\Meas{\Zahlmod{3}}$.  Let $\nu := \Phi\mu \in
\Meas{\Sof}$, so that if $\gU \subset \Sof$ is measurable, then $\nu [ \gU]
:= \mu \lb[ {\Phi}^{-1} (\gU)\rb]$.

 For all $N\in\Natur$, define character $\bchi_{N}$ by  $\D
\bchi_{N}(\bsft) := \prod_{n=0}^{N} (-1)^{x_{n}}$ for all $\bx\in\lb(\Zahlmod{2}\rb)^\Zahl$. Then Lemma \ref{even.odd} implies:
\beq
 \inn{ \bchi_{N}, \nu} 
& = & \nu (\Even_{N}) - \nu (\Odd_{N}) \\
& = &
\mu \lb( \gX_{0,0}^{N} \disj \gX_{1,2}^{N} \disj \gX_{2,1}^{N} \disj \gX_{0,2}^{N} \disj \gX_{1,0}^{N} \rb) 
\ - \ \nu \lb( \gX_{1,1}^{N} \disj \gX_{0,1}^{N} \disj \gX_{2,0}^{N} \disj \gX_{2,2}^{N}\rb).
\eeq
But $\mu$ is mixing, so
$\D \lim_{N\goto\oo} \mu (\gX_{i,j}^{N}) \ = \ \rho_{i}\cdot \rho_{j}$.
Thus, \
$\D \lim_{N\goto\oo} \inn{ \bchi_{N}, \nu} 
 \ = \ \rho_{0}^{2} + 2\rho_{1}\rho_{2} - \rho_{1}^{2} - \rho_{2}^{2}$. \  
So for example if 
 \[ \bP = \Matrix{ 1/2 & 0 & 1/2 \\
		   1/2 & 0 & 1/2 \\
		   0   & 1 & 0 \\}\]
with Perron measure $\brho \ = \ 
\lb(\frac{2}{5},\frac{1}{5},\frac{2}{5}\rb)$, then
$\D \lim_{N\goto\oo} \inn{ \bchi_{N}, \nu} \neq 0$. But clearly,
$\rank{\bchi_N}=N$, so that $\D \lim_{N\goto\oo} \rank{\bchi_N} \ = \
\oo$.  Thus $\nu$ is not harmonically mixing.

\section{Dispersion Mixing
\label{S:wkhm}}

  The example from \S\ref{S:non.harm.mixing} suggests the need for an
asymptotic randomization condition on measures that is less
restrictive than harmonic mixing.  In this section, we'll define the
concepts of {\em dispersion mixing} (for measures) and {\em
dispersion} (for automata) which together yield asymptotic
randomization.  In \S\ref{S:disperse} we'll show that many LCA are
dispersive.  In
\S\ref{S:unimix} and \S\ref{S:mkwrd} we'll show that many measures
 (including the Even Shift measure $\nu$ from
\S\ref{S:non.harm.mixing}) are dispersion mixing.

 Let $\Phi$ be an LCA as in equation (\ref{LCA}).
The advantage of this `polynomial' notation  is that
composition of two  LCA corresponds to multiplication of their
respective polynomials.
  For example, suppose  $\sA=(\Zahlmod{p})^s$, where $p\in\Natur$ is  prime, and $s\in\Natur$. Suppose $\Lat=\Zahl$ and 
$\Phi = 1+\shift{}$; \ that is, $\Phi(\ba)_0 = a_0 + a_1
\pmod{p}$.  Then the Binomial Theorem implies:
\beqn
\label{ledrappier.binom}
\mbox{For any $N\in\Natur$,}\qquad
\Phi^N \quad  = \quad  \sum_{n=0}^N \lb[N \atop n\rb]_p \shift{n},
\qquad\mbox{where}\quad 
\lb[N \atop n\rb]_p \ := \ \lb(N \atop n\rb) \bmod{p}.
\eeqn

  Let $S>0$, and let $\dK,\dJ\subset\Lat$ be subsets.  We say that
$\dK$ and $\dJ$ are {\dfn $S$-separated} if 
\[
\min\set{|\fk-\fj|}{\fk\in\dK \And \fj\in\dJ}
\quad \geq \quad S
\]
 If $\dF,\dG\subset\Lat$, and $\Phi \ = \ \D \sum_{\ff\in\dF}
 \varphi_\ff\cdot \shift{\ff}$ and $\Gam \ = \ \D \sum_{\fg\in\dG}
 \gam_\fg \cdot \shift{\fg}$ are two LCA, then we say $\Phi$ and
 $\Gam$ are {\dfn $S$-separated} if $\dF$ and $\dG$ are $S$-separated.
 Likewise, if $\dK,\dX\subset\Lat$, and $\bchi \ = \ \D
 \Tensor_{\fk\in\dK} \chi_\fk$ and $\bxi \ = \ \D \Tensor_{\fx\in\dX}
 \xi_\fx$ are two characters, then we say $\bchi$ and $\bxi$ are {\dfn
 $S$-separated} if $\dK$ and $\dX$ are $S$-separated.

 If $\Phi \ = \ \D \sum_{\ff\in\dF} \varphi_\ff\cdot \shift{\ff}$ is an LCA,
then let $\Srank{\Phi}$ be the maximum number of $S$-separated LCA
which can be summed to yield $\Phi$.  That is:
\[
 \Srank{\Phi} \ \ := \ \  \max\set{R}{\maketall \exists\, \Phi_1,\ldots,\Phi_R
\ \mbox{ mutually $S$-separated, with} \ \Phi \ = \ \Phi_1 + \cdots + \Phi_R}.
\]
For example, if 
\[ \Phi \quad = \quad 
1 \ + \ \shift{5}+\shift{6} \ + \ \shift{11}+\shift{12}+\shift{13},
\]
 then $\Srank[4]{\Phi} = 3$, because
$\Phi = \Phi_1 + \Phi_2 + \Phi_3$, where 
\[
 \Phi_1 \ = \ 1, \quad \Phi_2 \ = \ \shift{5}+\shift{6}, \And \Phi_3 \ = \
\shift{11}+\shift{12}+\shift{13}.
\]
  On the other hand, clearly, $\Srank[1]{\Phi} = 6$, while  
$\Srank[7]{\Phi} = 1$.

  Likewise, if $\bchi  \ = \ \D\Tensor_{\fk\in\dK} \chi_\fk$ is a
character, and $S>0$, then we define
\[
 \Srank{\chi}
\quad:=\quad \max\set{R}{\maketall \exists\, \bchi_1,\ldots,\bchi_R
\ \mbox{ mutually $S$-separated, with} \
 \bchi = \bchi_1 \tensor \cdots \tensor \bchi_R}.
\]
  (In the notation of \S\ref{S:mrf}, $\rank{\bchi} \ = \ \Srank[1]{\bchi}$.)
 
 We say that $\mu$ is {\dfn dispersion mixing} (DM) if,
for every $\eps>0$, there exist $S,R>0$ such that, for any
character $\bchi\in\h{\sA^\Lat}$, 
\quad
 $\statement{$\Srank[S]{\bchi} \ > \ R$}
\implies
\statement{$\lb|\inn{\bchi,\mu}\rb| \ < \ \eps$}$.
Note that dispersion mixing is less restrictive than harmonic mixing.

If $\Phi$ is an LCA and $\bchi$ is a character, then $\bchi\circ\Phi$
is also a character.
We say that $\Phi$ is {\dfn dispersive} if, 
for any $S>0$, and any character $\bchi\in\h{\sA^\Lat}$, 
there is a subset $\dJ\subset\dN$ of density 1 such that
$\D  \lim_{\dJ\ni j\goto \oo} \Srank{\bchi\circ\Phi^j} \  = \ \oo$.
It follows:

\Theorem{\label{AR.if.wkhm}}
{
  Let $\sA$ be any finite abelian group.
 If $\Phi:\sA^\dM\into\sA^\dM$ is a dispersive LCA and $\mu\in\sM(\sA^\dM)$ is dispersion mixing, then $\Phi$ asymptotically randomizes $\mu$.\qed
}
Theorem \ref{AR.if.wkhm} is an immediate
 consequence of an easily verified lemma:

\Lemma{\label{harmonic.AR.lemma}}
{
$\Phi$ asymptotically randomizes $\mu$
 if and only if,
for all $\bchi\in\h{\sA^\dM}$, there is a subset
$\dJ\subset\Natur$ with $\density{\dJ}=1$, such that
$\D  \lim_{\dJ\ni j\goto\oo}
\lb|\inn{\bchi\circ\Phi^{j}, \ \mu}\maketall\rb|
\quad=\quad 0$.
}
\bthmprf See the proof of Theorem 12 in \cite{PivatoYassawi1}. \ethmprf

\section{Dispersion and Bipartite CA\label{S:disperse}}

  If $\fm = (m_1,m_2,\ldots,m_D)\in\dM$, then let $|\fm| := |m_1| + 
|m_2| + \cdots + |m_D|$.  If $\Gam \ = \ \D \sum_{\fg\in\dG} \gam_\fg \cdot
\shift{\fg}$ is a linear cellular automaton, then define 
$\diam{\Gam} \ := \ \max\set{|\fg-\fh|}{\fg,\fh\in\dG}$.

 The {\dfn centre} of $\Gam$ is the centroid of $\dG$
 (as a subset of $\Real^n$):
\[
  \centre{\Gam} \quad:=\quad \frac{1}{\card{\dG}}\, \sum_{\fg\in\dG} \fg.
\]
  We say $\Gam$ is {\em centred} if $\lb|\centre{\Gam}\rb|<1$.
For any prime $p\in\Natur$, let
\[
  K_p\quad  := \quad \min\lb\{ \frac{1}{2}, \ \frac{4p-7}{4p+4}\rb\}.
\quad
\mbox{Thus,} \ \ K_2 \ = \ \frac{1}{12},\quad
\  K_3 \ = \ \frac{5}{16},\quad\mbox{and}\quad
   K_p \ = \ \frac{1}{2}, \ \ \mbox{for $p\geq 5$}.
\]
 Let  $\sA:=(\Zahlmod{p})^s$ (where $p$ is prime and $s\in\Natur$). 
 If $\Phi:\sA^\dM\into\sA^\dM$ is an LCA, then we say $\Phi$ is {\dfn bipartite} if $\Phi \
 = \ 1 + \Gam\circ
\shift{\ff}$, where $\Gam$ is centred and
 $\diam{\Gam} \ \leq \  K_p \cdot |\ff|$.
 For example:
\[
\begin{array}{rclcll}
  \Phi &=& 1+\shift{\ff} & 
\multicolumn{3}{l}{\mbox{is bipartite for any nonzero $\ff\in\dM$ and any prime $p\in\Natur$.}}\\
  \Phi &=& 1+\shift{12} + \shift{13} &= &
	 1 + \lb(1+\shift{}\rb)\circ \shift{12}&
\mbox{is bipartite for any prime $p\in\Natur$.}\\
  \Phi &=& 1+\shift{14} + \shift{19} &= & 
	 1 + \lb(\shift{-2}+\shift{3}\rb)\circ \shift{16}&
\mbox{is bipartite for any prime $p\geq 3$.}\\
  \Phi &=& 1+\shift{2} + \shift{3} &= & 1 + \lb(1+\shift{}\rb)\circ \shift{2}&
\mbox{is bipartite for any prime $p\geq 5$.}
\end{array}
\]
  Our goal in this section is to prove:

\Theorem{\label{bipart.disperse}}
{
Let $\sA=(\Zahlmod{p})^s$, where $p$ prime and $s\in\Natur$. 
  If $\Phi$ is bipartite, then $\Phi$ is dispersive.\qed
}

\breath

  For any $N\in\Natur$, let $\mtrx{N^{(i)}}{i=0}{\oo}$
 denote the {\dfn $p$-ary expansion} of $N$, so that
$\D  N \ =\  \sum_{i=0}^\oo N^{(i)} p^i$.
 Let $\Lucas{N} \ := \  \set{n\in\CC{0...N}}{n^{(i)} \leq
N^{(i)}, \mbox{for all} \ i \in \Natur}$.

\Lemma{(Lucas's Theorem) \label{lucas1}}
{
\bthmlist
  \item $\D \lb[N \atop n\rb]_p \ = \
\prod_{i=0}^\oo
  \lb[N^{(i)} \atop n^{(i)}\rb]_p$, \ 
where we define $\lb[N^{(i)} \atop n^{(i)}\rb]_p:=0$ if $n^{(i)}>N^{(i)}$,
and $\lb[0 \atop 0 \rb]_p:=1$.

   \item  Thus, $\lb[ N \atop n \rb]_p \not= \ 0$ iff $n \in \Lucas{N}$. \qed
\ethmlist
}
  For example, suppose $\Lat=\Zahl$ and 
$\Phi = 1+\shift{}$.  If we interpret equation (\ref{ledrappier.binom})
in the light of  Lemma \ref{lucas1}, we get:
\quad $\D \Phi^N \  = \  \sum_{n\in\Lucas{N}} \lb[N \atop n\rb]_p \shift{n}$.

\begin{figure}
\centerline{\includegraphics[angle=-90,scale=0.65]{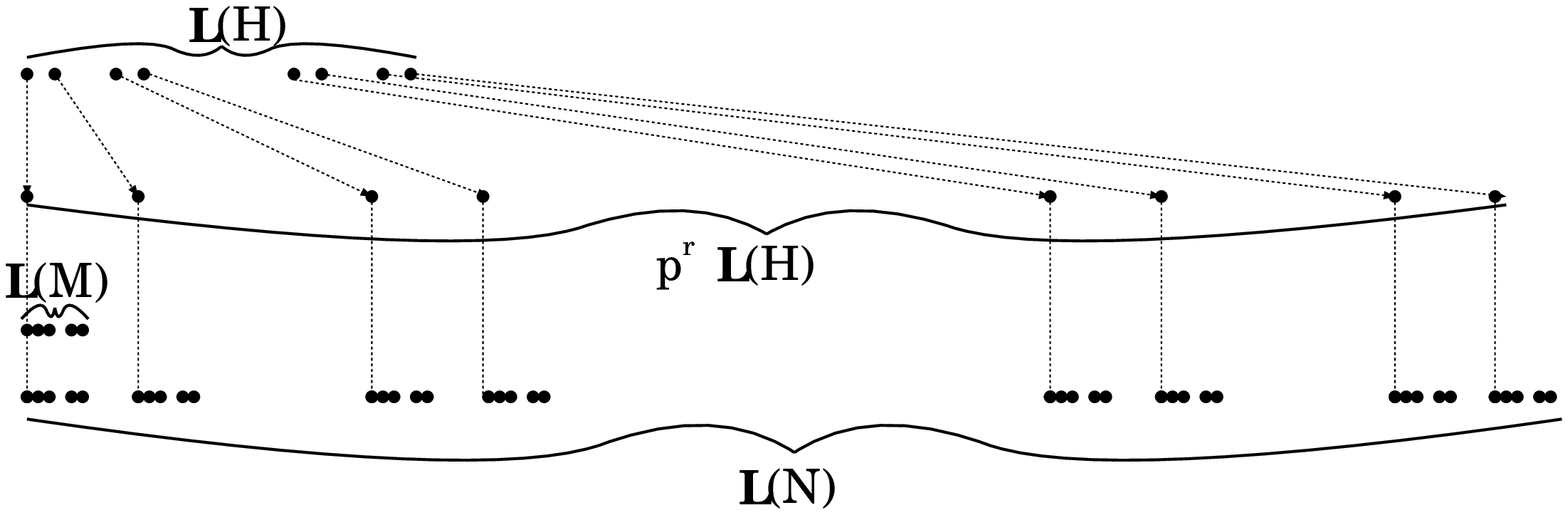}}
\caption{Lemma \ref{lucas.split}. \label{fig:lucas}}
\end{figure}

\Lemma{\label{lucas.split}}
{
Let $r, \ H \ \in \ \Natur$. 
\bthmlist
\item  If $M< p^r$, and $N=M + p^r\cdot H$, then \ 
$\D \Lucas{N} 
\ = \  \Lucas{M} + p^r\cdot\Lucas{H}$ (see Figure \ref{fig:lucas}).

\item If $m\in\Lucas{M}$, \  $h\in\Lucas{H}$,
and $n = m + p^r\cdot h$, then \ 
$\D
\lb[N\atop n\rb]_p  \ = \  
\lb[M\atop m\rb]_p \cdot \lb[H\atop h\rb]_p$.
\qed
\ethmlist
}

For example, suppose $p=2$ and $N \ = \ 53 \  =  \ 5 + 48 \ = \ 
5 +  2^4\cdot 3$.  Then $M=5$, $r=4$, and $H=3$, and
\beq
\Lucas{53} &=&
 \Lucas{5} + 2^4 \cdot \Lucas{3}
\quad=\quad \{0,1, \ 4,5\} \ \ + \ \  16\cdot\{0,1,2,3\}
\\&=&
\{ 0,1,\ 4,5, \quad 16,17, \ 20,21, \quad 32,33, \ 37,38, \quad 
	48,49, \ 52,53 \}.
\eeq
If $\bchi \ = \ \D\Tensor_{\fk\in\dK} \chi_\fk$ is a character, then
define $\diam{\bchi} \ := \ \max\set{|\fk-\fj|}{\fk,\fj\in\dK}$.
It follows:

\Lemma{\label{S.rank.lemma}}
{
  Let $\Phi$ be an LCA, and let $S>0$.
\bthmlist
\item If $\bchi$ is a character, and  $S_0 = S+\diam{\bchi}$, then 
 $\Srank{\bchi\circ \Phi} \ \geq \ \Srank[S_0]{\Phi}$.

\item If $\Gam$ is an LCA, and $S_0 = S+\diam{\Gam}$, then
 $\Srank{\Gam\circ \Phi} \ \geq \ \Srank[S_0]{\Phi}$.
\qed
\ethmlist
} 

\Corollary{\label{S.rank.lemma.2}}
{
 $\Phi$ is dispersive if and only if, for any $S_0>0$,
there is a subset $\dJ\subset\dN$ of density 1 such that
\ $\D 
  \lim_{\dJ\ni j\goto \oo} \Srank[S_0]{\Phi^j} \  = \ \oo$.\qed
}

 To prove Theorem \ref{bipart.disperse}, we'll use Lemma
\ref{lucas.split} to verify the condition of Corollary
\ref{S.rank.lemma.2}.  For any $S_0>0$, define
\[
 \dJ(S_0)  \quad := \quad \set{N\in\Natur}{N \ = \ M_N + p^{r_N} H_N, \ 
\mbox{for some} \ H_N,r_N > 0 \ \mbox{such that} \  M_N, \ S_0 < p^{r_N-1}}.
\]
For example, if $p=2$ and $S_0=7$, then $53 \in \dJ(7)$, because
$53 =  5 +  2^4\cdot 3$, so that $M_{53}=5$, $r_{53}=4$, and $H_{53}=3$.  Thus,
$2^{r_{53}-1} = 2^3=8$, and $7 < 8$ and $5<8$.
Note that $53 \ = \ 2^0 + 2^2 + 2^4 + 2^5$; thus, 
$53^{(3)}=0$.  This is exactly why $53\in\dJ(7)$:

\Lemma{\label{dJ.S0}}
{$\dJ(S_0) \  = \
 \set{N\in\Natur}{N \geq p\cdot S_0, \And N^{(r)} = 0 \ \mbox{for some} \ 
r\in \OC{\log_p(S_0)\ldots\log_p(N)}}$.}
\bthmprf
Suppose $N \ = \ M_N + p^{r_N} H_N$, 
for some $H_N, r_N > 0$ and $M_N\geq 0$, such that $M_N, \ S_0 < p^{r_N-1}$.
Let $r := r_N-1$; \  then $N^{(r)} = 0$ and $\log_p(S_0)< r < \log_p(N)$.

  Conversely, suppose $N^{(r)} = 0$, where $\log_p(S_0) < r <
\log_p(N)$.  Let $r_N := r+1$; \ then $S_0 \ < \ p^r \ = \ p^{r_N-1}$.
Let $M_N \ := \ \D\sum_{i=0}^{r-1} N^{(i)} p^i$; then $M_N  \ < \ p^r \ = \
p^{r_N-1}$ also.  Now let $H_N \ := \ \D\sum_{i=r_N}^\oo N^{(i)} p^{i-r_N}$;
then $N \ = \ M_N + p^{r_N} H_N$.
\ethmprf

\Lemma{\label{dJ.S0.density.1}}
{ $\density{\dJ(S_0)} = 1$.}
\bthmprf  
Let $\dI :=  \CC{p S_0\ldots \oo}$.  Then $\dI$ is a set
of density one, and Lemma \ref{dJ.S0} implies that
\[
 \dI\setminus\dJ(S_0) \quad=\quad \set{N\in\dI}
{ N^{(r)} \neq 0 \ \mbox{for all} \ r \in  \OC{\log_p(S_0)\ldots\log_p(N)}},
\]
which is a set of density zero.  It follows that  
$\density{\dJ(S_0)} = \density{\dI} = 1$.
\ethmprf

\Lemma{\label{lucas.power.splitting}}
{   If $N\in\dJ(S_0)$, and $N = M + p^r H$, then
\ $\Phi^N \ = \ \Phi^M \circ \THET^H$, where $\THET = \Phi^{(p^r)}$.
}
\bthmprf  Recall that $\Phi \ = \ 1 + \Gam\circ\shift{\ff}$.  Thus,
\beq
\Phi^N &\eeequals{(L)} &
  \sum_{n\in\Lucas{N}} \lb[N \atop n\rb]_p 
\lb(\maketall \Gam\circ \shift{\ff }\rb)^n
\quad\eeequals{(\ddagger)} \quad
 \sum_{m\in\Lucas{M}} \ \sum_{h\in\Lucas{H}}  \lb[H \atop h\rb]_p  \lb[M \atop m\rb]_p \lb(\maketall \Gam\circ \shift{\ff }\rb)^{(m + p^r h)} 
\\ &= &
 \sum_{h\in\Lucas{H}}  \lb[H \atop h\rb]_p \lb(\sum_{m\in\Lucas{M}} \   \lb[M \atop m\rb]_p \lb(\Gam \circ \shift{\ff}\rb)^m \rb) \circ 
\lb(\maketall \Gam\circ \shift{\ff }\rb)^{h p^r}
\\ &\eeequals{(\dagger)}&
 \sum_{h\in\Lucas{H}}  \lb[H \atop h\rb]_p  \Phi^M \circ
\lb(\maketall \Gam\circ \shift{\ff }\rb)^{p^r h}
\quad \eeequals{(\star)}\quad
 \Phi^M \circ \THET^H.
\eeq
{\bf(L)} is by Lucas Theorem and
$(\ddagger)$ is by Lemma \ref{lucas.split}(b).
\ 
{$(\dagger)$} is because $\D \Phi^M \ = \ \sum_{m\in\Lucas{M}} \
\lb[M \atop m\rb]_p \lb(\Gam \circ \shift{\ff}\rb)^m$.
Finally, 
{$(\star)$}  is because $\THET 
\  = \ \lb(1+ \Gam\circ\shift{\ff}\rb)^{p^r} \ \eeequals{(L)} \ 
  1+ (\Gam\circ\shift{\ff})^{p^r}$. Thus,
$\D \THET^H \ \eeequals{(L)} \
 \sum_{h\in\Lucas{H}}  \lb[H \atop h\rb]_p \lb(\maketall \Gam\circ \shift{\ff }\rb)^{p^r h}$.
\ethmprf

\bthmprf[Proof of Theorem \ref{bipart.disperse}.]
It suffices to verify the condition of Corollary
\ref{S.rank.lemma.2}.  So, let $S_1:=S_0+\diam{\Phi^M}$.  Then 
\beqn
\label{bipart.disperse.E2}
\Srank[S_0]{\Phi^N}
\quad\eeequals{(*)}\quad
\Srank[S_0]{\Phi^M\circ\THET^H}
\quad \geeeq{(\dagger)} \quad 
\Srank[S_1]{\THET^H}.
\eeqn
where $(*)$ is by Lemma \ref{lucas.power.splitting}
and $(\dagger)$ is by Lemma \ref{S.rank.lemma}(b).
    Thus, we want to show that
$\Srank[S_1]{\THET^H}\goesto{H\goto\oo}{} \oo$ for $H$ in a set of density 1.
To do this, we'll use {\em gaps} in $\Lucas{H}$.
If $h_0,h_1\in\Lucas{H}$, we say that $h_0$ and $h_1$ {\dfn bracket a  gap} if:
\[
 {\bf(i)} \ \  h_1 \geq p \cdot h_0
\qquad\And\qquad {\bf(ii)} \ \
\CO{h_0...h_1}\intsct\Lucas{H} \ = \ \emptyset.
\]
\begin{figure}
\centerline{
\includegraphics[angle=-90,scale=0.53]{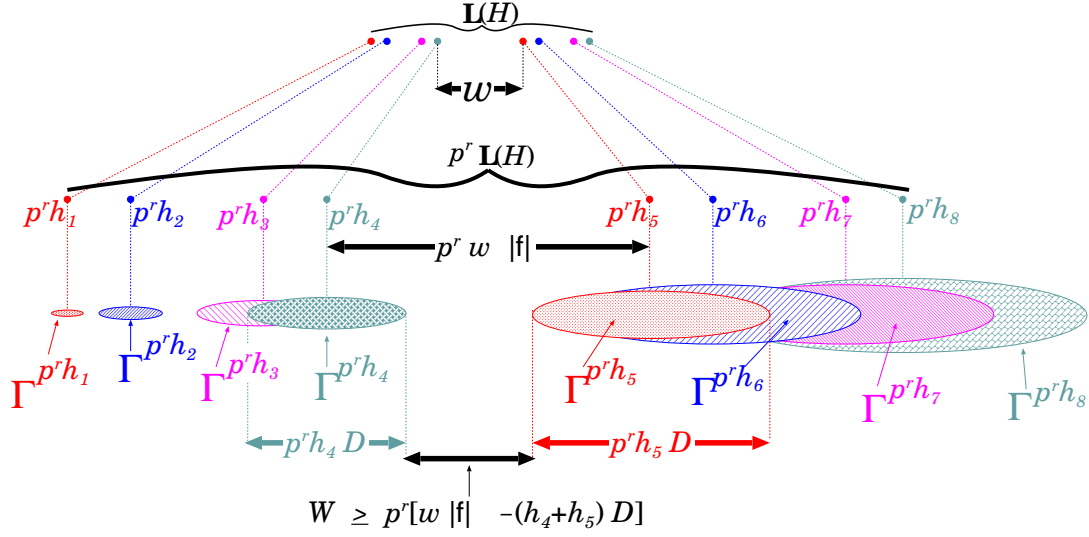}}
\caption{Claim \ref{bipart.disperse.C1} \label{fig:gaps}
of Theorem \ref{bipart.disperse}.}
\end{figure}

\Claim{\label{bipart.disperse.C1}
Let $h_0,h_1\in\Lucas{H}$, with $p\leq h_0<h_1$, and suppose
$h_0$ and $h_1$  bracket a gap in $\Lucas{H}$.
Then $\lb(\maketall \Gam\circ \shift{\ff }\rb)^{p^r h_0}$
and $\lb(\maketall \Gam\circ \shift{\ff }\rb)^{p^r h_1}$ are
$S_1$-separated. 
\ignore{Then $\THET^H$ has a gap greater than
$S_1$ between $\lb(\maketall \Gam\circ \shift{\ff }\rb)^{p^r h_0}$
and $\lb(\maketall \Gam\circ \shift{\ff }\rb)^{p^r h_1}$.}}
\bclaimprf
Suppose $|h_0-h_1| \ = \ w$. 
Then $\lb(\shift{\ff}\rb)^{p^r h_0}$ and $\lb(\shift{\ff}\rb)^{p^r h_1}$.
are  $(p^r\cdot w \cdot |\ff|)$-separated.
Thus, if $D=\diam{\Gam}$, then 
$\lb(\maketall \Gam\circ \shift{\ff }\rb)^{p^r h_0}$
and $\lb(\maketall \Gam\circ \shift{\ff }\rb)^{p^r h_1}$ are
$W$-separated, where
\begin{eqnarray}
\nonumber W&:=&
p^{r} w |\ff| \ -  \ \lb(\diam{ \Gam^{p_r h_0}}  \ +  \ \diam{ \Gam^{p_r h_1}}\rb)
\quad=\quad
p^{r} w |\ff| \ -  \  \lb(p^r h_0 D \ +  \  p^r h_1 D\rb)
\\ \label{E1} &\geq&
p^r\cdot\lb(w |\ff| - D\cdot (h_1+h_0)\maketall \rb).
\end{eqnarray}
(see Figure \ref{fig:gaps}).
We want $W \ \geq \ S_1$, or, equivalently, $W-\diam{\Phi^M} \ \geq \ S_0$
(because $S_1=S_0+\diam{\Phi^M}$).  First, note that
\begin{eqnarray}
\nonumber
\diam{\Phi^M}&\leq& M\cdot |\ff| \ + \ 2\cdot\max_{m\in\Lucas{M}}
\diam{\Gam^m} \quad= \quad M\cdot |\ff| \ + \ 2M\cdot D
\\&=&
 M\cdot \lb(|\ff| \ + \ 2D\maketall\rb)
\label{E2}
\quad\leq\quad
 p^{r-1} \cdot \lb(|\ff| \ + \ 2D\maketall\rb).
\end{eqnarray}
Thus,
\beq
W-\diam{\Phi^M}
&\geeeq{(*)} &
p^r\cdot\lb(w \cdot|\ff| - D\cdot (h_1+h_0)\maketall \rb)
\ - \ 
 p^{r-1} \cdot \lb(|\ff| \ + \ 2D\maketall\rb)
\\&=&
p^{r-1}\cdot\lb( p w \cdot|\ff| \ - \ p D\cdot (h_1+h_0)
\ - \  |\ff| \ - \ 2D\maketall\rb)
\\ &\geeeq{(\dagger)}&
S_0\cdot\lb( p w \cdot|\ff| \ - \ p D\cdot (h_1+h_0)
\ - \  |\ff| \ - \ 2D\maketall\rb).
\eeq
where $(*)$ is by equations (\ref{E1}) and (\ref{E2}), and
$(\dagger)$ is because $S_0<p^{r-1}$.

Thus, it suffices to show that
\[
   p w \cdot|\ff| \ - \ p D\cdot (h_1+h_0) \ - \  |\ff| \ - \ 2D \quad \geq \quad 1.
\]
To see this, observe that
\beq
\lefteqn{  p w \cdot|\ff| \ - \ p D\cdot (h_1+h_0) \ - \  |\ff| \ - \ 2D}\\
&=&
\lb( p w - 1\rb) \cdot |\ff| - \lb[p\cdot (h_1+h_0) - 2\maketall \rb]\cdot D
\quad\geeeq{(\flat)}\quad 
\lb( p w - 1\rb) \cdot |\ff| -
 \lb[p\cdot (h_1+h_0) - 2\maketall \rb]\cdot K_p\cdot|\ff|
\\&=&
\lb(\maketall p w - 1 - \lb[p\cdot (h_1+h_0) - 2\rb] K_p \rb) \cdot |\ff|
\quad\geeeq{(*)}\quad
 p\cdot(h_1-h_0)  - 1 - \lb[p\cdot (h_1+h_0) - 2\rb] K_p
\\&=&
 p\cdot\lb(\maketall (1-K_p)\cdot h_1 - (1+K_p)\cdot h_0 \rb) \ - \ (1+2\cdot K_p)
\\ & \geeeq{(\dagger)} & 
 p\cdot\lb(\maketall (1-K_p)\cdot p - (1+K_p) \rb)\cdot h_0 \ - \ 2
\quad \geeeq{(\ddagger)} \quad
 p^2\cdot\lb(\maketall (1-K_p)\cdot p - (1+K_p) \rb) \ - \ 2
\\ & \geeeq{(\star)} &
 \frac{3}{4}p^2 \ - \ 2
\quad \geeeq{(\diamond)} \quad 3 - 2 \quad=\quad 1.
\eeq
\bdesc
\item[$(\flat)$] is by hypothesis that $\Gam$ is bipartite.
\qquad{$(*)$} \ is because $|\ff|\geq 1$, and $w = h_1-h_0$.

\item[$(\dagger)$] is because $h_1 \geq p\cdot h_0$, and $K_p\leq \frac{1}{2}$.
\qquad{$(\ddagger)$}\ is because $h_0 \geq p$.
\item[$(\star)$] is because
$K_p \ \leq \ \frac{4p-7}{4p+4} \ = \ \frac{p-\frac{7}{4}}{p+1}$,  \ 
thus, $(p+1) K_p \ \leq \  p-\frac{7}{4} \ = \ p - 1 - \frac{3}{4}$; \ 
thus, $ \frac{3}{4} \ \leq \ (p - 1)  - (p+1) K_p 
\ = \ (1-K_p) p - (1+K_p)$.

\item[$(\diamond)$] is because $p\geq 2$, so $p^2\geq 4$.
\edesc
It follows that $W-\diam{\Phi^M} \geq S_0$, so that $W\geq S_1$.
\eclaimprf

  Let $\rank{H} \ := \ $ \# of gaps in $\Lucas{H}$.  Then 
Claim \ref{bipart.disperse.C1} implies that 
\beqn
\label{bipart.disperse.E3}
  \Srank[S_1]{\THET^H} \quad \geq \quad \rank{H}.
\eeqn
Thus, we want to show that the number of gaps is large.

  Suppose $i<k$.
We say that $i$ and $k$ {\dfn bracket a zero-block}
in the  $p$-ary expansion of $H$ if $H^{(i-1)}\neq 0 \neq H^{(k)}$,
but $H^{(j)}=0$, for all $i\leq j < k$.  For example, suppose $p=2$
and $H=19$.  Then $3$ and $5$ bracket a zero block in
the binary expansion $...010011$.

\Claim{\label{bipart.disperse.C2}
If $i$ and $k$ bracket a zero-block in the $p$-ary expansion of $H$,
then $p^i$ and $p^j$ bracket a gap in  $\Lucas{H}$.
}
\bclaimprf
 $H^{(i)}=0$, so the largest element in $\Lucas{H}$ less than
$p^i$ is 
\[
  h_0\quad=\quad
\sum_{j=1}^{i-1} H^{(j)}\cdot p^j
\quad\leq\quad
\sum_{j=1}^{i-1} (p-1)\cdot p^j \quad = \quad p^i - 1.
\]
 Now, $k \ = \ \min\set{j>i}{H^{(j)}\neq 0}$, so
 $h_1=p^k$ is the smallest element in $\Lucas{H}$ greater than
 $p^i$.  Also, $h_1 \ \geq \  p^{i+1} \  >\  p\cdot (p^i-1) \ \geq \
  p\cdot h_0$.  
\eclaimprf
\newcommand{\ZB}[1]{\#\bZ\bB\lb(#1\rb)}
Let $\ZB{H} \ := \   \# \mbox{of zero-blocks in the $p$-ary expansion of $H$}$.

Then Claim \ref{bipart.disperse.C2}  implies that 
\beqn
\label{bipart.disperse.E4}
\rank{H} \quad \geq \quad \ZB{H}.
\eeqn 
 Define \ $\dH \ := \ \set{H \in \Natur}{ \ZB{H} \ \geq \  \frac{1}{p^3} \log_p(H)}$.

\Claim{\label{density.H.is.one}
$\density{\dH}=1$.}
\bclaimprf
 Observe that $\ZB{H}$ is no less than the number of occurrences of
the word ``$101$'' in the $p$-ary expansion of $H$ (because $101$ is
a zero-block).  Let 
\[
\dH' \quad:=\quad \set{H \in \Natur}{
(\mbox{\# of occurrences of ``$101$''})  \ \geq \  \frac{1}{p^3} \log_p(H)}.
\]
Then $\dH'\subset\dH$. The Weak Law of Large Numbers implies
$\density{\dH'} = 1$.
\eclaimprf
 Define $\dJ  := \ \set{N\in\dJ(S_0)}{N = M_N + p^{r_N} H_N, \ \mbox{where} \
r_N \leq \frac{1}{2} \log_p(N), \ \mbox{and} \ H_N\in\dH}$.

\Claim{$\density{\dJ}=1$.}
\bclaimprf
$\dJ = \dJ_1\intsct\dJ_2$, where
\beq
\dJ_1 &:=& \set{N\in\dJ(S_0)}{N = M_N + p^{r_N} H_N, \ \mbox{where} \ 
 H_N\in\dH}\\
\And 
\dJ_2 &:=& \set{N\in\dJ(S_0)}{N = M_N + p^{r_N} H_N, \ \mbox{where} \
r_N \leq \frac{1}{2} \log_p(N)}.
\eeq
Now, $\density{\dJ_1}=1$ by Lemma \ref{dJ.S0.density.1} and 
Claim \ref{density.H.is.one}.  To see that $\density{\dJ_2}=1$, note that
\[
 \dJ(S_0)\setminus \dJ_2
\quad\subset\quad 
\set{N\in\Natur}{N^{(r)}\neq 0 \ 
\mbox{for all} \ r \in \OC{\log_p(S_0)\ldots\frac{1}{2} \log_p(N)}},
\]
  which is a set of density zero.
\eclaimprf
If $N = M_N + p^{r_N} H_N$ is an element of $\dJ$, then 
\beqn
\label{bipart.disperse.E5}
  \log_p(H_N) \quad\geq\quad \log_p(N) - r_N
\quad\geq\quad
 \log_p(N) - \frac{1}{2}\log_p(N) 
\quad=\quad \frac{1}{2}\log_p(N). 
\eeqn
Thus,
\beq
\Srank[S_0]{\Phi^N}
& \geeeq{(\heartsuit)} &
\Srank[S_1]{\THET^{H_N}}
\quad \geeeq{(\diamondsuit)}\quad
 \rank{H_N} \quad \geeeq{(\clubsuit)} \quad
 \ZB{H_N}  
\\&\geeeq{(*)}&
\frac{1}{p^3} \log_p(H_N) \quad \geeeq{(\spadesuit)} \quad
  \frac{1}{2p^3} \log_p(N).
\eeq
Here, $(\heartsuit)$ is by equation (\ref{bipart.disperse.E2}),\
$(\diamondsuit)$ is by equation (\ref{bipart.disperse.E3}), \  
$(\clubsuit)$ is by equation (\ref{bipart.disperse.E4}), \ 
$(\spadesuit)$ is by equation (\ref{bipart.disperse.E5}), 
and $(*)$ is because $H\in\dH$ by hypothesis.

 Thus $\D \lim_{\dJ\ni N\goto\oo}\, \Srank[S_0]{\Phi^N}
\quad \geq \quad \frac{1}{2p^3} \ \lim_{\dJ\ni N\goto\oo}\, \log_p(N)
\ = \ \oo$.
\ethmprf

\section{Uniform Mixing and Dispersion Mixing
\label{S:unimix}}

A measure $\mu\in\sM(\sA^\Zahl)$ is {\dfn uniformly mixing} if, for any
$\eps>0$, there is some $M>0$ such that, for any cylinder subsets
$\gL\subset\sA^\OC{-\oo...0}$ and $\gR\subset\sA^\CO{0...\oo}$, and any
$m>M$,
\beqn
\label{unimix.eqn}
  \mu\lb[\shift{m}(\gL) \ \intsct \gR\rb]
\quad\closeto{\eps}\quad
 \mu\lb[\gL\rb]\cdot\mu\lb[\gR\rb]
\eeqn
(here, ``$x \ \closeto{\eps} \ y$'' means $|x-y|<\eps$.)
\examples{ \label{X:unimix}
\item  Any mixing $N$-step Markov chain is uniformly mixing.
(See \S\ref{S:mkwrd}).

\item   If $\nu\in\sM(\sB^\Zahl)$ is uniformly mixing, and $\Psi:\sB^\Zahl\into \sA^\Zahl$ is a block map, then $\mu:=\Phi(\nu)$ is also uniformly mixing.
(If $\Psi$ has local map $\psi:\sB^\CC{-\ell...r}\into\sA$, then
replace the $M$ in (\ref{unimix.eqn}) with $M+\ell+r+1$).

\item Hence, if $\SFT\subset\sB^\Zahl$ is an SFT, and $\Sof:=\Psi(\SFT)\subset
\sA^\Zahl$ a sofic shift, and $\nu\in\sM(\SFT)$ is any mixing
$N$-step Markov chain, then $\mu:=\Phi(\nu)$ is a uniformly mixing 
measure on $\Sof$.  We call
$\mu$ a {\dfn quasi-Markov measure}.  } We say that $\mu$ is {\dfn
harmonically bounded} (HB) if there is some $C<1$ such that
$\lb|\inn{\bchi,\mu}\rb|<C$ for all $\bchi\in\h{\sA^\Zahl}$ except
$\bchi=\chr{}$.  The goal of this section is to prove:

\Theorem{
\label{unimix.wkhm}}
{
 Let $\sA$ be a finite abelian group.
 If $\mu\in\sM(\sA^\Zahl)$ is uniformly mixing and harmonically bounded, then $\mu$ is dispersion mixing.\qed
}
  We will then apply Theorem \ref{unimix.wkhm} to get:

\Corollary{\label{sofic.AR}}
{
 Let $\sA=\Zahlmod{p}$, where $p$ is prime.  If
$\mu\in\sM(\sA^\Zahl)$ is a mixing quasi-Markov measure, then
$\mu$ is asymptotically randomized by any dispersive LCA.  \qed
}

\paragraph*{Harmonic boundedness and entropy:}

\Lemma{\label{entropy.HB}}
{
  Let $\sA=(\Zahlmod{p})^s$, where $p$ is prime and $s\in\Natur$.
  If $\mu\in\sM(\sA^\Zahl)$ and
 $h(\mu,\shift{})\ > \ (s-1)\cdot\log_2(p)$, then $\mu$ is harmonically bounded.
}
\bthmprf
  Suppose $\mu$ was not HB.  Then for any $\alp>0$, we can find
$\chr{}\neq\bchi\in\h{\sA^\Zahl}$ with $|\inn{\bchi,\mu}|\ > \ 1-\alp$.
Let $\sI:=\image{\bchi}\subset\Torus{1}$, and let
$\nu:=\bchi(\mu)\in\sM(\sI)$ be the projected measure on $\sI$.  Thus,
$\inn{\bchi,\mu} \ = \ \D\sum_{i\in\sI} i\cdot \nu\{i\}$.
The following four claims are easy to check.

\Claim{\label{entropy.HB.claim.1}
For any $\bet>0$, there exists $\alp>0$ such that, for any
probability measure $\nu\in\sM(\sI)$ with
$\lb| \D\sum_{i\in\sI} i\cdot \nu\{i\}\rb|>1-\alp$, there
is some $i_0\in\sI$ with $\nu\{i_0\}>1-\bet$.\ecp}

\quad Suppose $\bchi=\D\Tensor_{k\in\dK}\chi_k$, where
$\dK\subset\CC{0...K}$ and $K\in\dK$.  Thus, if 
$\bxi :=\D \Tensor_{k\in\dK\setminus\{K\}}\chi_k$, 
then $\bchi = \bxi\tensor\chi_K$.  For any
$\bb\in\sA^\CO{0...K}$, let $\mu^{(\bb)}_K$ be the conditional measure
on the $K$th coordinate, and let $\nu^{(\bb)}_K \ := \
\chi_K\lb(\mu^{(\bb)}_K\rb)\in\sM(\sI)$ be the projected measure
on $\sI$.

\Claim{\label{entropy.HB.claim.2}
For any $\gam>0$, there exists $\bet>0$ such that,
if $\ \exists \, i_0\in\sI$ with $\nu\{i_0\}>1-\bet$, then
there is  a subset $\gB\subset\sA^\CO{0...K}$ with $\mu[\gB]>1-\gam$,
such that, for every $\bb\in\gB$,  there is some $i_\bb\in\sI$ with 
$\nu^{(\bb)}_K\{i_\bb\} \ > \ 1-\gam$. \
  Thus, if $\sP_\bb = \chi_K^{-1}\{i_\bb\}\subset\sA$, then
$\mu^{(\bb)}_K[\sP_\bb] \ > \ 1-\gam$.

(Observe that $\card{\sP_\bb}\leq p^{s-1}$
 for all $\bb\in\sA^\CO{0...K}$.) \ecp}

\quad For any measure $\rho\in\sM(\sA)$, define $\D   H(\rho) \ := \
 -\sum_{a\in\sA} \rho\{a\} \log_2\lb(\maketall \rho\{a\}\rb)$.
  Recall (e.g. \cite[Proposition 5.2.12]{Petersen})  that the $\shift{}$-entropy of $\mu$ can be computed:
\beqn
\label{entropy.limit}
  h(\mu,\shift{}) \quad = \quad \lim_{N\goto\oo} \ \sum_{\bb\in\sA^\CO{0...N}} \
\mu\cyl{\bb}\cdot H\lb(\mu^{(\bb)}_N\rb)
\eeqn

\Claim{\label{entropy.HB.claim.3}
For any $\del>0$, there exists $\gam_1>0$ such that, for any
probability measure $\rho$ on $\sA$, if there is a subset
$\sP\subset\sA$ with  $\card{\sP}\leq p^{s-1}$ and
$\rho[\sP] \ > \ 1-\gam_1$, then 
$H(\rho) \ < \ (s-1)\cdot\log_2(p) \ + \ \del$.\ecp}

\Claim{\label{entropy.HB.claim.4}
For any $\eps>0$, and $S>0$,
there exist $\del,\gam_2>0$ such that, for any $K\in\Natur$ and
probability measure $\mu$ on $\sA^\CC{0...K}$, if there is a subset
$\gB\subset\sA^\CO{0...K}$ with $\mu[\gB]>1-\gam_2$, such that,
for all $\bb\in\gB$, \ $H\lb(\mu^{(\bb)}_K\rb)< S-\del$, then
$\D \sum_{\bb\in\sA^\CO{0...K}}
\mu\cyl{\bb}\cdot H\lb(\mu^{(\bb)}_K\rb) \ < \ S-\eps$.
\ecp}

\quad  Now, set $S:=(s-1)\cdot\log_2(p)$. For any $\eps>0$, find $\del,\gam_2>0$
as in Claim \ref{entropy.HB.claim.4}.  Then find $\gam_1>0$ as in
Claim \ref{entropy.HB.claim.3}, and let $\gam:=\min\{\gam_1,\gam_2\}$.
Next, find $\bet$ as in Claim \ref{entropy.HB.claim.2} and then find
$\alp$ as in Claim \ref{entropy.HB.claim.1}.  Finally, find
$\bchi\in\h{\sA^\Zahl}$ with $|\inn{\bchi,\mu}|\ > \ 1-\alp$.  It then
follows from Claims \ref{entropy.HB.claim.1}-\ref{entropy.HB.claim.4}
that $\D \sum_{\bb\in\sA^\CO{0...K}}
\mu\cyl{\bb}\cdot H\lb(\mu^{(\bb)}_N\rb) \ < \ (s-1)\cdot\log_2(p)-\eps$.
But the limit in (\ref{entropy.limit}) is a decreasing limit,
so we conclude that $h(\mu,\shift{}) \ < \ (s-1)\cdot\log_2(p)-\eps$.
Since this is true for any $\eps>0$, we conclude that 
$h(\mu,\shift{}) \ \leq \ (s-1)\cdot\log_2(p)$, contradicting our hypothesis.
\ethmprf

\Corollary{\label{zero.entropy.HB}}
{
  If $\sA=\Zahlmod{p}$ {\rm (where $p$ is prime)},
and $h(\mu,\shift{})>0$, then $\mu$ is harmonically bounded.\qed
}

Say $\mu$ is {\dfn uniformly multiply mixing} if, for any $\eps>0$,
there is some $S>0$ such that, for any $R>0$, if $\dK_0,\dK_1,\ldots,\dK_R\subset\Lat$ are finite, mutually $S$-separated subsets of $\Lat$, and
$\gU_0\subset\sA^{\dK_0},\ldots,\gU_R\subset\sA^{\dK_R}$ are
cylinder sets, then:
\beqn
\label{unimultimix.eqn}
  \mu\lb(\Intsct_{r=0}^R \gU_r\rb)
\quad\closeto{R\eps}\quad
 \prod_{r=0}^R\mu\lb(\gU_r\rb).
\eeqn

\Lemma{}
{
  If $\mu\in\sM(\sA^\Zahl)$
 is uniformly mixing, then $\mu$ is uniformly multiply mixing.
}
\bthmprf  ({\em by induction on $R$}).  
The case $R=1$ is just uniform mixing.  Suppose
(\ref{unimultimix.eqn}) is true for all $R'<R$.  Find $S>0$ so
that, if $\dK_0,\ldots,\dK_R$ are mutually $S$-separated, then
\[
 \mu\lb(\Intsct_{r=0}^R \gU_r \rb)
\quad = \quad  \mu\lb(\gU_0\intsct \ \Intsct_{r=1}^R \gU_r \rb) 
 \quad \closeto{\eps} \quad
\mu\lb(\gU_0\rb)\cdot  \mu\lb(\Intsct_{r=1}^{R} \gU_r \rb)
\quad \closeto{(R-1)\eps} \quad 
 \mu\lb(\gU_0\rb)\cdot \prod_{r=1}^R\mu\lb(\gU_r\rb),
\]
where ``$\closeto{\eps}$'' comes by setting  $R'=1$, and
``$\closeto{(R-1)\eps}$'' comes by setting $R'=R-1$.
\ethmprf
  
\Lemma{\label{function.mixing}}
{
 Suppose $\mu\in\sM(\sA^\Zahl)$ is uniformly multiply mixing.  For any $\eps>0$
and $R\in\Natur$, there is some $S>0$ such that:
if  $\dK_0,\ldots,\dK_R\subset \Zahl$ are  $S$-separated sets,
and, for all $r\in\CC{0...R}$, \quad
  $\bchi_r:\sA^{\dK_r}\into\Cplx$ are characters,
and  \  $\bchi \ = \ \D \prod_{r=0}^{R} \bchi_r$, \
then \ 
$\D
\inn{\bchi, \ \mu }
\quad\closeto{\eps/2}\quad
\prod_{r=0}^{R} \inn{\bchi_r, \ \mu }$.\qed
}

\bthmprf[Proof of Theorem \ref{unimix.wkhm}.]
Let $\eps>0$.  We want to find $S>0$ and $R>0$ such that, if
$\bchi$ is any character, and $\Srank[S]{\bchi}>R$, then $\lb|\inn{\bchi,\mu}\rb|<\eps$.

 Let $C<1$ be the  harmonic bound.
Find $R\in\Natur$ such that $C^R<\eps/2$.
Let $S>0$ be as in Lemma \ref{function.mixing}.
Suppose $\Srank[S]{\bchi}>R$, and
let $\bchi \ := \ \D \Tensor_{r=0}^R \bchi_r$, where  $\bchi_r:\sA^{\dK_r}\into\Cplx$ are characters, and $\dK_0,\ldots,\dK_R\subset\Zahl$ are $S$-separated.
Then Lemma \ref{function.mixing} implies:
\beqn
\label{mark.word.cond.exp}
\inn{\bchi, \ \mu }
\quad\closeto{\eps/2}\quad
\prod_{r=0}^{R} \inn{\bchi_r, \ \mu }.
\eeqn
By harmonic boundedness, we know
$ \lb|\inn{\bchi_r, \ \mu }\rb| < C$ for all $r\in\CC{0...R}$.
Thus, (\ref{mark.word.cond.exp}) implies:
$\D
\lb|\inn{\bchi, \ \mu }\rb|
\quad\closeto{\eps/2}\quad
\prod_{r=0}^{R} \lb|\inn{\bchi_r, \ \mu }\rb|
\quad < \quad
\prod_{r=0}^{R} C
\quad=\quad
C^{R+1}
\quad<\quad
 C^R
\quad<\quad \eps/2$.
\ethmprf

\bprfof{Corollary \ref{sofic.AR}}
  From examples \ref{X:unimix}(a) and \ref{X:unimix}(b), we know $\mu$ is
uniformly mixing.  Any mixing quasi-Markov measure has nonzero entropy,
so Corollary \ref{zero.entropy.HB} says that $\mu$ is
harmonically bounded.   Theorem \ref{unimix.wkhm}  says
$\mu$ is dispersion mixing.  Theorem \ref{AR.if.wkhm} says
$\mu$ is asymptotically randomized  by any dispersive CA.
\ethmprf
\section{Markov Words
\label{S:mkwrd}}

 If $m,n\in\Zahl$, and $m\leq n$, let $\sA^\CO{m...n}$ be the set of all
{\dfn words} of the form $\ba=[a_m, \ a_{m+1},\ldots,a_{n-1}]$.
Let $\D \sA^* := \Union_{-\oo<m< n<\oo} \sA^\CO{m...n}$ be the set of all finite
words.  Elements of $\sA^*$ are denoted by boldfaced letters (e.g.
$\ba$, $\bb$, $\bc$), and subsets by gothic letters (e.g. $\gA$, $\gB$,
$\gC$).  Concatenation of words is indicated by juxtaposition.
Thus, if $\ba=[a_0 \ldots a_n]$ and $\bb=[b_0 \ldots b_m]$, then
$\ba\bb = [a_0 \ldots a_n b_0 \ldots b_m]$.

If $V>0$ and  $\bv\in\sA^\CO{-V...V}$, we say that
$\bv$ is a {\dfn Markov word} for $\mu$ if (in the terminology of
\S\ref{S:mrf}), $\bv$ isolates $\OO{-\oo...-\!V}$
from $\CO{V...\oo}$.

\examples{
\item If
$\mu$ is an $N$-step Markov shift, and $N\leq 2V$, then every
$\bv\in\sA^\CO{-V...V}$ is a Markov word.

\item Let $\SFT\subset\sB^\Zahl$ be a subshift of finite type,
let $\Psi:\SFT\into\sA^\Zahl$ be a block map, so that $\Sof:=\Psi(\SFT)$
is a sofic shift.  Let $\nu$ be a Markov measure on $\SFT$ and
let $\mu:=\Psi(\nu)$.  If $\sof\in\Sof_\CC{-V...V}$ is a synchronizing
word for $\Psi$, then $\sof$ is a Markov word for $\mu$.
}

\Proposition{\label{prop.mkword.unimix}}
{  If $\mu\in\sM(\sA^\Zahl)$ is mixing and has a Markov word, then $\mu$ is
uniformly mixing.
}
\bthmprf  Fix $\eps>0$.
 For any words $\ba,\bb\in\sA^*$, the mixing of $\mu$ implies that
there is some $M_\eps\lb({\ba},{\bb}\rb)<\oo$ such that, for all
$m>M_\eps({\ba},{\bb})$,
\quad
  $\mu\lb(\shift{m}\cyl{\ba} \ \intsct\cyl{\bb}\rb)
\ \closeto{\eps} \  \mu\cyl{\ba}\cdot\mu\cyl{\bb}$.
   Our goal is to find some $M>0$ such that $M_\eps({\ba},{\bb})<M$ for
all $\ba,\bb\in\sA^*$.

Let $\bv\in\sA^*$ be a Markov word for $\mu$.

\Claim{\label{mkword.unimix.claim2}
Let $\bu,\bw,\bu',\bw'\in\sA^*$, and consider the words
$\bu\bv\bw$ and $\bu'\bv\bw'$.  We have: \ 
$ M_\eps\lb(\bu\bv\bw, \ \bu'\bv\bw'\rb)
\ = \ 
  M_\eps\lb(\bv\bw, \ \bu'\bv\rb)$.
}
\bclaimprf 
Define transition probabilities:
$  \mu(\bu \comesfrom \bv) \ := \ \mu(\bu\bv)/\mu(\bv)$
and
  $\mu(\bv \goesinto \bw) \ := \ \mu(\bv\bw)/\mu(\bv)$.
If $m> M_\eps\lb(\bv\bw, \ \bu'\bv\rb)$, then
\begin{eqnarray}
\mu\lb(\shift{m}\cyl{\bu\bv\bw}\intsct\cyl{\bu'\bv\bw'}\maketall \rb)
\label{mkword.unimix.claim2.1}
&= &
\mu(\bu \comesfrom \bv)\cdot
\mu\lb(\shift{m}\cyl{\bv\bw}\intsct\cyl{\bu'\bv}\maketall \rb)\cdot
\mu(\bv \goesinto \bw')\qquad\\
\label{mkword.unimix.claim2.2}
&\closeto{\eps}&
\mu(\bu \comesfrom \bv)\cdot
\mu\cyl{\bv\bw}\cdot\mu\cyl{\bu'\bv}\cdot
\mu(\bv \goesinto \bw')\\
\label{mkword.unimix.claim2.3}
&=&
\mu\cyl{\bu\bv\bw}\cdot\mu\cyl{\bu'\bv\bw'}.
\end{eqnarray}
  (\ref{mkword.unimix.claim2.1}) and (\ref{mkword.unimix.claim2.3})
 are because $\bv$ is a Markov word; \ 
(\ref{mkword.unimix.claim2.2}) is because $m> M_\eps\lb(\bv\bw,  \bu'\bv\rb)$.
\eclaimprf

  If $\ba\in\sA^*$, we say that $\bv$ {\dfn occurs}
in $\ba$  if $\ba\restr{\CO{n\!-\!V...n\!+\!V}} \ = \ \bv$
for some $n$.  

\Claim{\label{mkword.unimix.claim1}
There is some $N>0$ such that
$\mu\set{\ba\in\sA^\CC{0...N}}{\mbox{$\bv$ occurs in $\ba$}}
\ > \ 1-\eps$.}
\bclaimprf
By ergodicity, find $N$ such that $\D\mu\lb(\Union_{n=0}^N
\shift{n}\cyl{\bv}\rb) \ > \ 1-\eps$.
\eclaimprf
Let $\sA^*_\bv$ be the set of words (of length at least $N$)
in $\sA^*$ with $\bv$ occuring in
the last $(N+V)$ coordinates, and let $_\bv\sA^*$ be the set of all words
in $\sA^*$ with $\bv$ occuring in the first $(N+V)$ coordinates.
 Then Claim \ref{mkword.unimix.claim1} implies that:
\beqn
\label{mkword.unimix.eqn2}
\mu(\sA^*_\bv)\quad>\quad 1-\eps \quad \And \quad
\mu(_\bv\sA^*)\quad>\quad 1-\eps.
\eeqn 
 Let $\D\sA^{<N} \ := \  \Union_{n=1}^N \sA^\CC{0...n}$.  Then
\beqn
\label{mkword.unimix.eqn1}
\begin{array}{rcl}
\sA^*_\bv
&=&
\set{\bu\bv\bw}{  \bu\in\sA^* \And \bw\in\sA^{<N}}.
\\  \And 
_\bv\sA^*
&=&
\set{\bu'\bv\bw'}{ \bu'\in\sA^{<N} \And \bw'\in\sA^*},
\end{array}
\eeqn
\beq
\mbox{ Define}\quad  M_1 &:=&  \max_{\ba \in \sA^*_\bv}\ \max_{\bb\in _\bv\sA^*} \
M_\eps(\ba,\bb)
\quad\eeequals{(*)}\quad
\max_{{\scriptstyle \bu\in\sA^*} \atop {\scriptstyle \bw\in\sA^{<N}}} \
\max_{{\scriptstyle \bu'\in\sA^{<N}}\atop{\scriptstyle \bw'\in\sA^*} } 
 M_\eps\lb(\bu\bv\bw, \ \bu'\bv\bw'\rb)
\\&\eeequals{(\dagger)}&
\max_{\bw,\bu'\in\sA^{<N}} \ M_\eps\lb(\bv\bw, \ \bu'\bv\rb).
\eeq
 where $(*)$ is by equation  (\ref{mkword.unimix.eqn1})
and $(\dagger)$ is by Claim \ref{mkword.unimix.claim2}.
Likewise, define
\beq
 M_2 &:=& \max_{\ba \in \sA^*_\bv}\ \max_{\bb\in \sA^{<N}} \ M_\eps(\ba,\bb)
\quad=\quad \max_{\bw \in \sA^{<N}}\ \max_{\bb\in \sA^{<N}} \
 M_\eps(\bv\bw, \ \bb),\\
 M_3 &:=& \max_{\ba \in \sA^{<N}}\ \max_{\bb\in _\bv\sA^{*}} \ M_\eps(\ba,\bb)
\quad=\quad  \max_{\ba\in \sA^{<N}} \ \max_{\bu' \in \sA^{<N}}\
 M_\eps(\ba, \ \bu'\bv),\\
\And M_4 &:=& \max_{\ba \in \sA^{<N}}\ \max_{\bb\in \sA^{<N}} \ M_\eps(\ba,\bb).
\eeq
Thus, $M_1,\ldots,M_4$ each maximizes a finite collection of finite
values, so each is finite.  Thus, $M:=\max\{M_1,\ldots,M_4\}$ is finite.

\Claim{For any $\ba,\bb\in\sA^*$, \quad $M_\eps(\ba,\bb) <M$.}
\bclaimprf
If $\ba \in \sA^{<N} \union  \sA^*_\bv$ and 
 $\bb \in \sA^{<N} \union  \  _\bv\sA^*$, then $M_\eps(\ba,\bb) <M$ by
definition.

So, suppose $\ba \not\in \sA^{<N} \union  \sA^*_\bv$.  Then 
equation  (\ref{mkword.unimix.eqn2}) implies that $\mu[\ba]<\eps$.  Hence,
for any $m\in\Natur$,\quad $\mu(\shift{m}[\ba]\intsct\bb) \ < \  \eps$
and
 $\mu[\ba]\cdot \mu[\bb] \ < \  \eps$.
Thus, $ \mu(\shift{m}[\ba]\intsct\bb) \closeto{\eps} \mu[\ba]\cdot \mu[\bb]$ automatically.  Hence, $M_\eps(\ba,\bb) = 0  <M$.

  Likewise, if  $\bb \not\in \sA^{<N} \union  \  _\bv\sA^*$, then 
 $M_\eps(\ba,\bb) = 0  <M$.
\eclaimprf
Thus, $\mu$ is uniformly mixing.
\ethmprf

\Corollary{\label{hb.mix.mrk.wd.random}}
{
  If $\mu$ is harmonically bounded, mixing and has a Markov word, 
then $\mu$ is asymptotically randomized by $\Phi=1+\shift{}$.
}
\bthmprf
Combine Proposition \ref{prop.mkword.unimix} with Theorems
\ref{AR.if.wkhm} and \ref{unimix.wkhm}. \ethmprf

\section{Lucas Mixing\label{S:lucas.mix}}

  Throughout this section, let $D:=1$, so that
$\Lat=\Zahl$.  Let $\sA:=(\Zahlmod{p})^s$, where $p\in\Natur$ is
prime, and $s\in\Natur$.  Let $\Phi := 1 + \shift{}$.
 We will introduce a condition on $\mu$ which is weaker
than dispersion mixing, and which is both sufficient and {\em
necessary} for asymptotic randomization.

  Let $\bchi\in\h{\sA^\Zahl}$, and suppose
$\bchi \ = \ \D \Tensor_{\fk\in\dK} \chi_\fk$.  We define
$ \dm{\bchi} \ := \ \max(\dK)-\min(\dK)$, and
define 
\[
 \ldm{\bchi} \quad := \quad p^r,
\quad \mbox{where}\quad r \ := \  \logdiam{\bchi}.
\]
It follows from Lucas' Theorem that $\Phi^\ldm{\bchi} \ = \ 1 + \shift{\ldm{\bchi}}$.  Thus, for any $h\in\Natur$, 
\[
\Phi^{h\cdot \ldm{\bchi}} \ = \quad  \sum_{\ell\in\Lucas{h}}
\lb[h\atop \ell\rb]_p \shift{\ldm{\bchi}\cdot\ell}, \quad
\mbox{and thus,}\quad
\dilate{\bchi}{h}  \ = \quad \Tensor_{\ell\in\Lucas{h}}
\lb[h\atop \ell\rb]_p \bchi\circ\shift{\ldm{\bchi}\cdot\ell}.
\]
 Observe that $\dK+p^r\ell$ and $\dK+p^r\ell'$ are
disjoint for any $\ell\neq\ell'\in\Lucas{h}$.   Hence, if $L:=\card{\Lucas{h}}$, then
$\dilate{\bchi}{h}$ is a product of $L$ `disjoint translates' of $\bchi$.

  If $\mu$ is a measure on $\sA^\Zahl$, we say that $\mu$ is 
{\dfn Lucas mixing} if, for any nontrivial character $\bchi\in\h{\sA^\Zahl}$,
there is a subset $\dH\subset\Natur$ of \Cesaro density one such that
\ $\D \lim_{\dH\ni h\goto\oo} \inn{\dilate{\bchi}{h}, \ \mu} \ = \  0$.
Our goal in this section is to prove:

\Theorem{\label{lucas.mix}}
{
 $\statement{$\Phi=1+\shift{}$ asymptotically randomizes $\mu$}
\iff
\statement{$\mu$ is Lucas mixing}$.\qed 
}

  It is relatively easy to see that:

\Lemma{}
{
  If $\mu$ is dispersion-mixing, then $\mu$ is Lucas mixing.\qed
}

  Thus, the ``$\seilpmi$'' direction of Theorem \ref{lucas.mix} is
an extension of Theorem \ref{AR.if.wkhm}, in the case 
$\Phi=1+\shift{}$.  The ``$\implies$'' direction 
makes this the strongest possible extension for this LCA.

\breath

  Set $S := \dm{\bchi}$, and let $\tldJ \ := \ \dJ(S)$, where
$\dJ(S)$ is defined as in \S\ref{S:disperse}.  It follows from
Lemma \ref{dJ.S0.density.1} that $\density{\tldJ}=1$.
  For any $m\in\Natur$, let $\bchi^m \ := \ \bchi\circ\Phi^m$.

\Lemma{\label{dilation.lemma}}
{
 Let $j\in\tldJ$, with $j \ = \ m + p^r \cdot h$.
Then $\bchi\circ\Phi^j \ = \ \dilate{\bchi^m}{h'}$,
where $h' = p^s \cdot h$ for some $s\geq0$.
}
\bthmprf
Apply Lemma \ref{lucas.power.splitting} to observe that 
$\Phi^j \ = \ \Phi^m \circ \Phi^{h\cdot (p^r)}$.
Thus, 
\[
\bchi\circ\Phi^j \quad = \quad \bchi\circ\Phi^m \circ \Phi^{h\cdot (p^r)}
\quad = \quad \bchi^m \circ \Phi^{h\cdot (p^r)}.
\]
By definition, $r$ is such that $m<p^{r-1}$ and
$\dm{\bchi}<p^{r-1}$.  Thus, 
\[
\dm{\bchi^m}  \quad =  \quad \dm{\bchi}+m
  \quad < \quad  p^{r-1}+p^{r-1} 
  \quad \leq \quad  p^r.
\]
  Now, let $s \ := \ r- \log_p\dm{\bchi^m}$, and let $h' \ := \  p^s \cdot h$.
 Then $h\cdot (p^r)\ = \ h'\cdot\ldm{\bchi^m}$, so that
$\Phi^{h\cdot (p^r)} \ = \ \Phi^{h'\cdot\ldm{\bchi^m}}$.
\ethmprf

\bthmprf[Proof of Theorem \ref{lucas.mix}.]
 We will use Lemma \ref{harmonic.AR.lemma}.

`$\seilpmi$'\quad  
For any $m\in\Natur$, let $r(m) \ := \
\lb\lceil\log_p\lb(\max\lb\{m, \ \dm{\bchi}\maketall\rb\}\rb)\rb\rceil
 \ + \ 1$,
and define
\beqn
\label{natur.m.defn}
\tldJ_m \quad := \quad  \set{ m + p^{r(m)} h}{ h\in\Natur}.
\eeqn
It follows that:
\beqn
\label{natur.bchi.union}
\tldJ \quad=\quad
\Union_{m\in\Natur} \tldJ_m.
\eeqn
If $j \ = \ m + p^{r(m)}h$ is an element of $\tldJ_m$, then
Lemma \ref{dilation.lemma} says $\bchi\circ\Phi^j \ = \ \dilate{\bchi^m}{h'}$,
for some $h'\geq h$.
 Now, $\mu$ is Lucas mixing, so find a subset $\tldH_m\subset\Natur$ of density one with
$\D \lim_{\tldH_m\ni h\goto\oo} \inn{\dilate{\bchi^m}{h}, \ \mu}
\ = \ 0$. Define:
\begin{eqnarray}
\nonumber
 \dH_m &:=&
 \set{h\in\tldH_m}{\lb|\inn{\dilate{\bchi^m}{h}, \ \mu}\rb|
 \ \leq \ \frac{1}{m}}, \\
\label{natur.prime.m.defn}
 \dJ_m &:=&
 \set{m + p^{r(m)} h}{h\in\dH_m}, \\
\And
\label{natur.prime.bchi.union}
 \dJ & := & \D \Union_{m\in\Natur} \dJ_m.
\end{eqnarray}
\Claim{ $\density{\dJ}=1$.}
\bclaimprf
For any $m\in\Natur$, there is some $K$ such that
$\dH_m  =  \tldH_m\intsct\CO{K...\oo}$.  Thus, $\reldense{\dH_m/\tldH_m}
 =  1$.  Thus, $\density{\dH_m}  =  \density{\tldH_m}  =  1$. 
Compare (\ref{natur.m.defn}) and (\ref{natur.prime.m.defn})
to see that  $\reldense{\dJ_m/\tldJ_m}  =  1$. 
Then compare (\ref{natur.bchi.union}) and (\ref{natur.prime.bchi.union})
to see that $\reldense{\dJ/\tldJ}  =  1$.
Thus, $\density{\dJ}  =  \density{\tldJ} =  1$. \eclaimprf

\Claim{$\D \  \lim_{\dJ\ni j\goto\oo} \inn{\bchi\circ\Phi^j, \ \mu}
\ = \ 0$.}
\bclaimprf
  Fix $\eps>0$.  Let $M$ be large enough that $\frac{1}{M}<\eps$.  For
all $m\in\Natur$ with $m<M$, find $H_m$ such that, if $h\in\tldH_m$ and
$h>H_m$, then $\lb|\inn{\dilate{\bchi^m}{h}, \ \mu}\rb| \ < \ \eps$.
Let $J_m:= m + 2^{r (m)} \cdot H_m$.  Thus, if
$j\ = \ m+2^{r (m)}\cdot h$ is an element of
$\dJ_m$, and $j \ > \ J_m$, then we must have $h>H_m$, so
that $\lb|\inn{\bchi\circ\Phi^j, \ \mu}\rb| \ = \ \lb|\inn{\dilate{\bchi^m}{h}, \ \mu}\rb| \ < \ \eps$.

\qquad Now let $J \ := \ \D \max_{1\leq m \leq M} \ J_m$.  
Thus, for all $j\in\dJ$, if $j>J$, then either 
$j\in\dJ_m$ for some $m\leq M$, in which case 
$\lb|\inn{\bchi\circ\Phi^j, \ \mu}\rb| \ < \ \eps$ by  construction of $J$,
or $j\in\dJ_m$ for some $m>M$, in which case
\[
\lb|\inn{\bchi\circ\Phi^j, \ \mu}\rb| \quad \lt{(*)} \quad
 \frac{1}{m} 
\quad<\quad \frac{1}{M} \quad  \lt{(\dagger)} \quad  \eps.
\]
Here, $(*)$ follows by definition of $\dH_m$, and $(\dagger)$ follows
by definition of $M$.
\eclaimprf
Lemma \ref{harmonic.AR.lemma} and Claims 1 and 2 imply
that $\Phi$ asymptotically randomizes $\mu$.

`$\implies$'\quad Suppose $\mu$ was not weakly harmonically mixing.
Thus, there is some $\bchi\in\h{\sA^\Zahl}$ and some subset
$\dH\subset\Natur$ of density $\del>0$ such that $\D \ \limsup_{\dH\ni
h\goto\oo} \ \lb|\inn{\dilate{\bchi}{h}, \ \mu}\rb|
\ > \ 0$.  But  $\dilate{\bchi}{h} \ = \ \bchi\circ\Phi^{p^r\cdot h}$ \ 
(where $r=\logdiam{\bchi}$).
Hence, if $\dJ \ := \ p^r\cdot\dH$, then
$\density{\dJ} = p^{-r}\cdot \del\ > \ 0$, and
$\D
\limsup_{\dJ\ni j\goto\oo}\ \lb|\inn{\bchi \circ \Phi^j, \ \mu}\rb|
\ = \ 
\limsup_{\dH\ni h\goto\oo} \ \lb|\inn{\dilate{\bchi}{h}, \ \mu}\rb|
\ > \ 0$.  But then Lemma \ref{harmonic.AR.lemma} implies that $\Phi$
cannot randomize $\mu$.
\ethmprf

\section{Randomization of Zero-Entropy Measures
\label{S:zero.random}}

   Of the probability measures which are
asymptotically randomized by LCA, every known
example has positive entropy.  However, we'll show that positive
entropy is {\em not} necessary, by constructing a class of zero-entropy
measures which are {Lucas mixing}, and thus (by Theorem \ref{lucas.mix})
randomized by $\Phi=1+\shift{}$.

  For both efficiency and lucidity, we will employ
probabilistic language.  Let
$(\Omg,\sB,\rho)$ be an abstract probability space (called the {\em
sample space}).  If $(\bX,\sX)$ is any measurable space, then an
($\bX$-valued) {\dfn random variable} is a measurable function
$f:\Omg\into\bX$.  In particular, a {\dfn random sequence} is a
measurable function $\ba:\Omg\into\sA^\Zahl$.  By
convention, we suppress the argument of random variables.  Thus, if
$\ba,\bb,\bc$ are random sequences, then the equation ``$\ba +\bb =
\bc$'' means ``$\ba(\omg) + \bb(\omg) = \bc(\omg)$, \ for
$\rho$-almost all $\omg\in\Omg$.'' 

 If $f:\Omg\into\bX$ is a random variable, 
and $\bU\subset\bX$,
then ``$\Prob{f\in\bU}$'' denotes $\rho\lb[f^{-1}(\bU)\rb]$.
If $g:\Omg\into\bY$ is another random variable, then 
$f$ and $g$ are {\dfn independent} if, for any measurable 
$\bU\subset\bX$ and $\bV\subset\bY$, 
\quad
$\Prob{f\in\bU \ \mbox{and} \ g\in\bV} \ = \ \Prob{f\in\bU}\cdot\Prob{g\in\bV}$
---i.e. $ \rho\lb[f^{-1}(\bU) \intsct g^{-1}(\bV)\rb]
\ = \ 
  \rho\lb[f^{-1}(\bU) \rb] \cdot  \rho\lb[ g^{-1}(\bV)\rb]$.
 The {\dfn distribution} of $f$ is the probability
measure $ \mu := f(\rho)$ on $(\bX,\sX)$; we then say that $f$ is a
{\dfn $\mu$-random variable}.  Thus, every random variable determines
a probability measure on its range.  However, given a measure $\mu$,
we can construct infinitely many independent $\mu$-random variables.

\breath

   Let $\sA:=\Zahlmod{2}$ and $\mu\in\sM(\sA^\Zahl)$, and consider 
a $\mu$-random sequence $\ba\in\sA^\Zahl$.
We say $\mu$ has {\dfn independent random dyadic increments} (IRDI) if,
for any  $n\in\Natur$, and all $m\in\CC{1...2^n}$,
 $a_{m+2^n} = a_{m} + d^n_m$, where $d^n_1,\ldots,d^n_{2^n}$ are independent
$\sA$-valued random variables.
If $d^n_1,\ldots,d^n_{2^n}$ have distributions $\del^n_1,\ldots,\del^n_{2^n}$,
then $\mu$ has {\dfn lower decay rate} $\alp\in\OO{0,1}$ if
there is some $L>0$ such that, for all $n\geq L$, and all $m\in\CC{1...2^n}$,
\quad $\alp^n \ \leq \ \del^n_m\{1\}$.

\Proposition{\label{IRDI.asympt.rand.lemma}}
{ 
If $\mu$ has IRDI with lower decay rate $\alp>\frac{1}{\sqrt{2}}$,
then $\mu$ is Lucas Mixing.
}
\bthmprf
Let $\bchi\in\h{\sA^\Zahl}$ be a nontrivial character.
We seek  $\dH\subset\Natur$ with $\density{\dH}=1$, such
that
$\D \lim_{\dH\ni h \goto\oo} \inn{\bchi\circ\Phi^{h\cdot \ldm{\bchi}}, \ \mu}
 \ = \ 0$.

  If $n\in\Natur$, let $I=I(n):=\lb\lceil\log_2(n)\rb\rceil$, and suppose $n$
has binary expansion
$\{n^{(i)}\}_{i=0}^{I}$.  Let $\dI(n) := \set{j\in\CC{0...I}}{n^{(j)} =1}$.
Let $\eps>0$ be small, and define:
\[
  \dH \quad:=\quad\set{h\in\Natur}{\card{\dI(h)} \  \ \geq \  \  \frac{1}{2} I(h) - \eps}.
\]
  Then $\density{\dH}=1$.  Suppose $n\in\dH$ is large;  let $\dI:=\dI(n)$
and $I:=I(n)$.  Assume $I$ is large (in particular, $I>L$).

Now, $\alp>\frac{1}{\sqrt{2}}$, so
find $\bet$ such that $\frac{1}{\alp} \ < \ \bet \ < \ \sqrt{2}$. 
Define
\beqn
\label{M.bigger.beta}
M\quad:=\quad \card{\dI}-1
 \quad\geq \quad
\frac{1}{2} I  - \eps - 1
\quad  \quad \grt{(*)} \quad \log_2(\beta) I,
\eeqn
where $(*)$ is  because $\log_2(\beta) \ < \ \frac{1}{2}$ and $I$ is large,
while $\eps$ is small.

  Suppose $\dI = \{ i_1 < i_2 < \ldots < i_{M+1} = I\}$. 
Let $\bxi_0:=\bchi$, and for 
each $m\in\CC{0...M}$, define \
 $\bxi_{m+1} :=  \bxi_m  \tensor \ \lb(\bxi_m\circ\shift{L_i}\rb)$, \
where $L_i:=2^{i_m} \cdot\ldm{\chi}$. \ 
Thus, \  $\bchi\circ \Phi^{n\cdot\ldm{\chi}}  \ =  \ \bxi_{M+1}$.  

  Let $r := \rank{\bchi}$.  Then for all $m\in\CC{1...M+1}$, \ \ 
$\rank{\bxi_m} \ = \ 2^m \cdot r$.  In particular, define
\beqn
 R \quad := \quad \rank{\bxi_M} \quad = \quad  2^M \cdot r  
\quad \grt{(*)} \quad
 \beta^I \cdot r. 
\label{R.defn}
\eeqn
where $(*)$ is by equation (\ref{M.bigger.beta}).
Thus, \
$\D  \bxi_M  \ =  \ \Tensor_{x\in\dX} \xi_x$, \ 
where $\dX\subset\Zahl$ is a subset with $\card{\dX}=R$.  Thus,
if $\ba\in\sA^\Zahl$ is a $\mu$-random sequence, then
\begin{eqnarray}\nonumber
\bxi_{M+1}(\ba) &=&
\bxi_{M}(\ba) \ \cdot \  \lb(\bxi_M\circ\shift{2^I}(\ba)\rb)
\quad = \quad  \prod_{x\in\dX} \xi_x(a_x) \cdot  \xi_x \lb(a_{x+2^I}\rb)
\\ &=&
\label{IRDI.asympt.rand.lemma.e1}
 \prod_{x\in\dX} \xi_x\lb(a_x + a_{x+2^I}\rb)
\quad = \quad \prod_{x\in\dX} \xi_x\lb(d^I_x\rb),
\end{eqnarray}
where $\{d^I_x\}_{x\in\dX}$ are independent random dyadic increments.
If $d^I_x$ has distribution $\del^I_x$, then
\beqn
 \Expct[\del^I_x]{ \xi_x\lb(d^I_x\rb)} 
\quad=\quad
\del^I_x\{0\}-\del^I_x\{1\} 
\quad=\quad 
1-2\del^I_x\{1\}  
 \quad \leeeq{(*)} \quad 
1-2\cdot\alp^{I}
\quad = \quad 
\frac{2\alp^{-I} - 1}{2\alp^{-I}}. 
\label{IRDI.asympt.rand.lemma.e2}
\eeqn
Here, $(*)$ is because $\mu$ has  lower decay rate $\alp$, so
$\del^I_x\{1\} \geq \alp^I$ (assuming $I\geq L$).
\[
\mbox{Thus,}\quad
\inn{\mu, \ \ \bchi\circ\Phi^n}
\quad\eeequals{(\ddagger)} \quad
\Expct[\mu]{\prod_{x\in\dX} \xi_x\lb(d^I_x\rb)}
\quad\eeequals{(*)} \quad
\prod_{x\in\dX} \Expct[\del^I_x]{ \xi_x\lb(d^I_x\rb)}
\quad\leeeq{(\dagger)} \quad
\lb(\frac{2\alp^{-I} - 1}{2\alp^{-I}} \rb)^R.
\]
Here, $(\ddagger)$ is by equation (\ref{IRDI.asympt.rand.lemma.e1}),
\  $(*)$ is because $\{d^I_x\}_{x\in\dX}$ are independent, and
$(\dagger)$ is by equation (\ref{IRDI.asympt.rand.lemma.e2}) and
because $\card{\dX}=R$.
\beq
\mbox{Thus,}\quad
\log \lb|\maketall \inn{\mu, \ \ \bchi\circ\Phi^n}\rb|
& \leq & 
 R \cdot \lb[\maketall 
\log\lb(2\alp^{-I} - 1\rb)-\log(2\alp^{-I})\rb] 
\quad \leeeq{(\ast)} \quad
 -R \cdot \log'\lb(2\alp^{-I}\rb)\\
& = & \frac{-R}{2\alp^{-I}}
\quad \lt{(\dagger)} \quad \frac{-\beta^I\, r}{2\alp^{-I}} 
\quad = \quad
- \frac{r}{2} \ (\alp\bet)^I.
\eeq
Here, $(\ast)$ is because $\log$ is a decreasing function,
and $(\dagger)$ is by equation (\ref{R.defn}).

But $\bet> \frac{1}{\alp}$, so  $\alp\bet>1$.
Thus,
$\D \lim_{\dH\ni h\goto\oo} \log  \lb|\maketall \inn{\mu, \  \bchi\circ\Phi^{h\cdot \ldm{\bchi}}}\rb|
\ = \ - \frac{r}{2} \ \lim_{I\goto\oo} \ 
(\alp\bet)^I \ = \ -\oo$.  Hence $\D \lim_{\dH\ni
h\goto\oo} \lb|\inn{\mu, \ \bchi\circ\Phi^{h\cdot \ldm{\bchi}}}\rb|
\  = \ 0$.
\ethmprf

Suppose $\mu\in\sM(\sA^\Zahl)$ has independent random dyadic increments;
for any  $n\in\Natur$, and all $m\in\CC{1...2^n}$,
let  $\del^n_1,\ldots,\del^n_{2^n}$ be
the dyadic increment distributions, as before.
Then $\mu$ has {\dfn upper decay rate} $\alp\in\OO{0,1}$ if
there are constants $L_1,K>0$ such that, for all
 $n\geq L_1$, and all $m\in\CC{1...2^n}$,
\quad $\del^n_m\{1\} \ \leq \ K\cdot \alp^{n}$.

\Proposition{\label{IRDI.entropy.lemma}}
{ 
  If $\mu$ has IRDI with upper decay rate $\alp<1$, then $h(\mu) = 0$.
}
\bthmprf
 Let $L_1, K>0$ be as above.  Assume without loss of generality that $K>4$.
Let $L_2 \ := \ \D\frac{-\log_2(K) - 1}{\log_2(\alp)}$. 
 Let $L:=\max\{L_1,L_2\}$.

For any $n\in\Natur$, and $m\in\CC{1...2^n}$,
let $\del^n_m$ be  as above.
The entropy of $\del^n_m$ is defined:
\beqn
\label{distropy.defn}
  H(\del^n_m) \quad  := \quad  -\del^n_m\{0\}\log_2(\del^n_m\{0\})  \ - \ \del^n_m\{1\}\log_2(\del^n_m\{1\})
\eeqn

\Claim{There exists  $c_1>0$ such that,
if $n>L$ and $m\in\CC{1...2^n}$, then
$H(\del^n_m) \  <  \  c_1 n\cdot\alp^{n}$.}
\bclaimprf
 $\alp<1$, so $\log_2(\alp)<0$;  Thus, if $n\geq L_2$, then
$n\log_2(\alp) \leq  L_2\log_2(\alp)$. Thus,
\begin{eqnarray}
\nonumber
\log_2(K\alp^n) &=&\log_2(K)+n\log_2(\alp)
\quad\leq \quad\log_2(K)+L_2 \log_2(\alp)
\\&=&
\label{IRDI.entropy.lemma.e2}
\log_2(K) - \log_2(K)-1 \quad= \quad-1.
\end{eqnarray}
Thus, 
$\del^n_m\{1\} \  \leeeq{(*)} \  K\alp^n  \  \leeeq{(\dagger)} \ \frac{1}{2}$,
where $(*)$ is because $n\geq L_1$ and $(\dagger)$ is
by equation (\ref{IRDI.entropy.lemma.e2}).

 But, 
if $\del^n_m\{1\}<\frac{1}{2}$ in equation (\ref{distropy.defn}), 
then $H(\del^n_m)$ decreases as $\del^n_m\{1\}$ decreases.
Hence, 
\beq
H(\del^n_m) & \leq & 
  -K\alp^{n}\log_2 \lb(K\alp^{n}\rb) \ - \
 \lb(1-K\alp^{n}\rb)\log_2 \lb(1-K\alp^{n}\rb)
\\& <&
  K \alp^{n}\underbrace{(nA-k)}_{(*)} \ + \ \lb(1-K\alp^{n}\rb)\cdot 
\underbrace{2 K\alp^{n}}_{(\dagger)}
 \quad = \quad
 K  \lb(nA+2-k-2K\alp^n\rb)\cdot \alp^{n}\\
 & \lt{(\ddagger)} &  KnA\cdot \alp^{n} 
\quad\lt{(\diamond)}\quad c_1 n\cdot \alp^n.
\eeq
Here, $(*)$ is the substitution $k:=\log_2(K)$ and $A:=-\log_2(\alp)$;
\quad
$(\dagger)$ is because, if $\eps$ is small, then $\log(1-\eps) \approx -\eps$,
thus,  $-\log(1-\eps) < 2\eps$;\quad
$(\ddagger)$ is because $2-k-2K\alp^n<0$ because $k>2$ because we assume
$K>4$;\quad
$(\diamond)$ is where $c_1:= KA>0$.
\eclaimprf

  Let $\ba\in\sA^\Zahl$ be a $\mu$-random sequence, and
fix $n>L$.
To compute the conditional entropy 
$H\lb(\ba\restr{\OC{2^n...2^{n+1}}} | \ba\restr{\CC{1...2^n}}\rb)$,
recall that, for all $m\in\CC{1...2^n}$,
\quad $a_{2^n+m} \ = \ a_m + d^n_m$.  Thus, 
\begin{eqnarray}
\nonumber
H\lb(\ba\restr{\OC{2^n...2^{n+1}}} \rb|\lb. \ba\restr{\CC{1...2^n}}\rb)
&=&
H\lb(d^n_1,d^n_2,\ldots,d^n_{2^n}\rb) 
\quad  \eeequals{(*)} \quad
\sum_{m=1}^{2^n} H(\del^n_m)
\\ & \lt{(\dagger)} &
 2^n\cdot c_1 n \alp^n
\quad  = \quad  c_1 n \cdot (2\alp)^n.
\label{IRDI.entropy.lemma.e3}
\end{eqnarray}
where $(*)$ is because  $d^n_1,d^n_2,\ldots,d^n_{2^n}$ are independent random variables with distributions $\del^n_1,\ldots,\del^n_{2^n}$,
and $(\dagger)$ is by Claim 1. 
\quad Thus, for any $N>L$,
\begin{eqnarray}\nonumber
H\lb(\ba\restr{\CC{1...2^N}} \rb|\lb. \ba\restr{\CC{1...2^{L}}} \rb)
&=&
\sum_{n=L}^{N-1} 
H\lb(\ba\restr{\OC{2^n...2^{n+1}}}  \rb|\lb. \ba\restr{\CC{1...2^n}}\rb) 
\quad \lt{(*)} \quad
\sum_{n=L}^{N-1}
   c_1 n \cdot (2\alp)^n \\
& < & \nonumber
 c_1 N \cdot (2\alp)^{L}\sum_{n=0}^{N-L-1}
  (2\alp)^n 
\quad = \quad
 c_1 N  \cdot (2\alp)^{L}\ \frac{(2\alp)^{N-L}-1}{2\alp -1}
\\ &\leq&
 c_2 N \cdot (2\alp)^N,
\label{IRDI.entropy.lemma.e4}
\end{eqnarray}
where  $(*)$ is by equation (\ref{IRDI.entropy.lemma.e3}), and
where  $c_2 \ \D \approx \ \frac{c_1}{2\alp-1} \ > \ 0$ is another constant.

Thus, if $H_0 :=  H\lb(\ba\restr{\CC{1...2^{L}}} \rb)$, then
\beqn
H\lb(\ba\restr{\CC{1...2^N}}  \rb)
\quad= \quad
H\lb(\ba\restr{\CC{1...2^N}} \rb|\lb. \ba\restr{\CC{1...2^{L}}} \rb) \ + \  H_0
\quad\leeeq{(*)} \quad
 c_2 N \cdot (2\alp)^N \ + \ H_0,\qquad\label{IRDI.entropy.lemma.e}
\eeqn
where $(*)$ is by equation (\ref{IRDI.entropy.lemma.e4}).  Thus,
\beq
 h(\mu) 
&=&
\lim_{M\goto\oo}\ \frac{1}{M}\ H\lb(\ba\restr{\CC{1...M}} \rb)
\quad= \quad
\lim_{N\goto\oo}\ \frac{1}{2^N}\ H\lb(\ba\restr{\CC{1...2^N}} \rb)\\
&\leeeq{(*)}&
\lim_{N\goto\oo}\ \frac{c_2 N \cdot (2\alp)^N + H_0 }{2^N} 
\quad\leq \quad 
  c_2 \lim_{N\goto\oo}\  N \alp^{N} 
\quad\eeequals{(\dagger)} \quad 0,
\eeq
where $(*)$ is by equation (\ref{IRDI.entropy.lemma.e}), and
$(\dagger)$ is because $|\alp|<1$.
\ethmprf

\begin{figure}
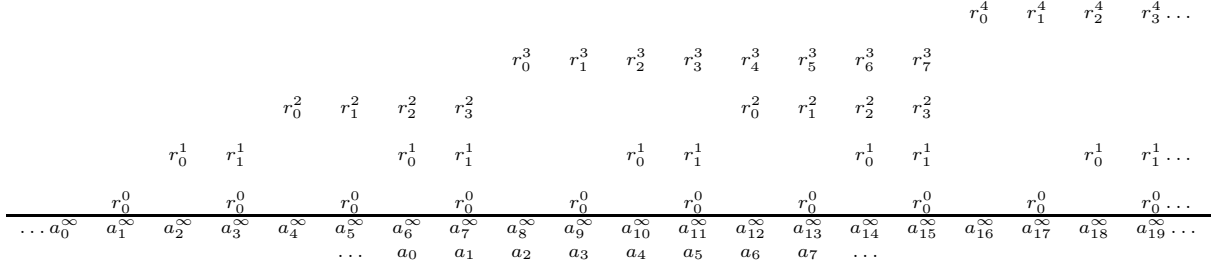

{\scriptsize
\[
\begin{array}{cccccccccccccccccccccccccccccccccccccccccccccccccccccc}
 & & & & & & & &  & & & & & & & & r^4_0 &r^4_1 &r^4_2 &r^4_3 \ldots\\\\
 & & & & & & & & r^3_0 &r^3_1 &r^3_2 &r^3_3 &r^3_4 &r^3_5 &r^3_6 &r^3_7 &
 & & &  \\ \\
 & & & &r^2_0&r^2_1&r^2_2&r^2_3& & & & &r^2_0&r^2_1&r^2_2&r^2_3& & & &\\\\ 
 & & r^1_0 & r^1_1 & & & r^1_0 & r^1_1 & & & r^1_0 & r^1_1 & & & r^1_0 & r^1_1 & & & r^1_0 & r^1_1 \ldots \\\\
&r^0_0& &r^0_0& &r^0_0& &r^0_0& &r^0_0& &r^0_0& &r^0_0& &r^0_0& &r^0_0& &r^0_0
\ldots\\
\hline
\ldots a^\oo_{0} & a^\oo_{1} & a^\oo_{2} & a^\oo_{3} & a^\oo_{4} & a^\oo_{5} & a^\oo_{6} & a^\oo_{7} & a^\oo_{8} & a^\oo_{9} & a^\oo_{10} & a^\oo_{11} & a^\oo_{12} & a^\oo_{13} & a^\oo_{14} & a^\oo_{15} & a^\oo_{16} & a^\oo_{17} & a^\oo_{18} & a^\oo_{19}  \ldots\\
&&&&&\ldots& a_{0} & a_{1} & a_{2} & a_{3} & a_{4} & a_{5} & a_{6} & a_{7} & \ldots
\end{array}
\]}
\caption{\footnotesize The construction of random sequence $\ba^{\oo}$;
the approximation of $\ba$ as a random translate of $\ba^{\oo}$. \label{fig:construct}}
\end{figure} 

  It remains to actually construct a measure with IRDI.
  Let $0 < \alp < 1 $.  For any $n\in\Natur$, let
$\rho_n$ be the probability distribution on $\sA=\Zahlmod{2}$ such that
\beqn
\label{prob.def}
  \rho_n\{1\} = \alp^{n} \qquad\And\qquad
 \rho_n\{0\} = 1-\alp^{n}.
\eeqn
  For each $n\in\Natur$, we will construct a random sequence $\ba^{n}
\in\sA^\Zahl$ as follows.  First, define $\ba^{0} :=
[\ldots0000\ldots]$.  Now, suppose, inductively, that we have
$\ba^{n}$.  Let $r^n_0,r^n_1,\ldots,r^n_{2^n-1}$ be a set of $2^n$ independent
$\sA$-valued,  $\rho_n$-random variables.  Let
$\br^{n}\in\sA^\Zahl$ be the random, $2^{n+1}$-periodic sequence
\[
  \br^{n} \ := \ [\ldots,
\underbrace{\stackrel{\lefteqn{\mathrm{\scriptstyle zeroth \ coordinate}}\atop\downarrow}{0},0,\ldots,0}_{2^n},
r^n_0,r^n_1,\ldots,r^n_{2^n-1},\underbrace{0,0,\ldots,0}_{2^n},r^n_0,r^n_1,\ldots,r^n_{2^n-1},\ldots],
\]
 and inductively define $\ba^{n+1} := \ba^{n} + \br^{n}$.

  Let $\mu_n\in\sM(\sA^\Zahl)$ be the distribution of $\ba^{n}$, and
let $\tl\mu_n := \D \frac{1}{2^n}\sum_{i=1}^{2^n} \shift{i}(\mu_n)$ be
the stationary average of $\mu_n$. Finally, let $\mu := \D
\wkstlim_{n\goto\oo} \tl\mu_n$.

  Let $\mu_{\oo}$ be the probability distribution of the
random sequence $\ba^{\oo} := \D \sum_{n=1}^\oo \br^{n}$ 
(see Fig.\ref{fig:construct}).  Then $\mu_\oo = \D \wkstlim_{n\goto\oo} \mu_n$, and
loosely speaking, $\mu$ is the `$\shift{}$-ergodic average' of $\mu_{\oo}$.
Thus, if $\ba$ is a $\mu$-random sequence, we can think of $\ba$ as
obtained by shifting $\ba^{\oo}$ by a random amount.  The next lemma
describes the structure of $\ba^\oo$:

\Lemma{\label{toplitz.lemma}}
{
  Let $M\in\Natur$ have binary expansion 
$\D M=\sum_{n=0}^\oo m_n 2^n$.  For all $n\geq0$, let $\D M_n:=
\sum_{i=0}^{n-1} m_i 2^i$.  Then 
$\D a^\oo_M \ = \ \sum_{n=0}^\oo m_n\cdot r^n_{M_n}$.\qed
}
  For example, suppose $M:=13 = 1 + 4 + 8$;  then $m_0=m_2=m_3=1$
and $m_1=0$.  Hence, $M_0=0$, \ $M_1=M_2 = 1$, and $M_3 = 5$.
Thus, $a^\oo_{13} \ = \ r^0_0 + r^2_1 + r^3_5$
(see Figure \ref{fig:construct}).

  Think of $\ba^{\oo}$  as being generated by a process
of `duplication with error'.  Let $\bw^{0} := [0]$ be a word of
length 1.  Suppose, inductively, that we have $\bw^{n} = [w_1
w_2\ldots w_{2^n-1}]$.  Let $\tl\bw^{n} := [\tlw_1
\tlw_2\ldots \tlw_{2^n-1}]$ be an `imperfect copy' of $\bw^{n}$: \quad
for each $m\in\CO{0...2^n}$, \ $\tlw_m := w_m + r^n_m$, where
$r^n_0,r^n_1\ldots,r^n_{2^n-1}$ are the independent $\rho_n$-distributed
variables from before, which act as `copying errors'.  Let $\bw^{n+1}
\ := \ \bw^n \tl\bw^n$.  Then
$\ba^{\oo}$ is the limit of $\bw^{n}$ as $n\goto\oo$.

\Proposition{\label{toplitz.irdi}}
{
  $\mu$ has IRDI, with upper and lower
decay rate $\alp$.
}
\bthmprf  
 Let $\ba\in\sA^\Zahl$ be a $\mu$-random sequence, and fix $N\in\Natur$.
By construction, there is some $k\in\Zahl$ such that
$\ba$ looks like $\shift{k}(\ba^{\oo})$ in a neighbourhood
around $0$.  To be precise, 
\beqn
\label{toplitz.irdi.e4}
\mbox{For all $m\in\CO{0...2^{N+1}}$,}\qquad  
a_{m}\quad =\quad a^{\oo}_{k+m}.
\eeqn
  For example, in Figure \ref{fig:construct}, let $N=2$,
so that $2^N=4$; suppose $k=6$.  Thus,
$ [a_0,a_1,\ldots,a_7] 
\ = \ 
 [a^\oo_6,a^\oo_7,\ldots,a^\oo_{13}]$.
\   Thus, $d^2_0 = a_4-a_0 \ = \ a^\oo_{10} - a^\oo_6 \ = \ 
r^3_2 - r^2_2  =  r^3_2 + r^2_2$.  More generally:

\Claim{\label{toplitz.irdi.c1} 
Let $m\in\CO{0...2^N}$. 
\bthmlist
\item There is a set  $S(m) := \{(n_0,m_0), \ (n_1,m_1),\ldots,(n_J,m_J)\}$
(for some $J\geq 0$), where  $N=n_0 \leq n_1 \leq \cdots \leq n_J$,
and where $m_j\in\CO{0...2^{n_j}}$ for $\forall \ j\in\CC{0...J}$,  such that 
\ $d^N_m   \ = \ 
 r^{n_0}_{m_0} +  r^{n_1}_{m_1} + \ldots +  r^{n_J}_{m_J}$.

\item If $m'\in\CO{0...2^N}$, and $m'\neq m$, then $S(m')\intsct S(m)=\emptyset$.
\ethmlist
}
\bclaimprf
Let $M:=k+m$ and let $\tlM:=k+m+2^N$.
If $\D M=\sum_{n=0}^\oo m_n 2^n$ and $\D \tlM=\sum_{n=0}^\oo \tlm_n 2^n$, 
then Lemma \ref{toplitz.lemma} says that
\beqn
\label{toplitz.irdi.e5}
a^\oo_M \quad = \quad \sum_{n=0}^\oo m_n\cdot r^n_{M_n};
\qquad\And\qquad
 a^\oo_{\tlM} \ = \ \sum_{n=0}^\oo \tlm_n\cdot r^n_{\tlM_n}.
\eeqn
Let  $N_1\geq N$ be the smallest element of $\CO{N...\oo}$ such that
 $m_{N_1}=0$.  Hence,  $m_n=1$ for all $n\in\CO{N...N_1}$, and $m_{N_1}=0$. 
Note that  $\tlM=M+2^N$, so binary expansions of $M$ and $\tlM$ are related as follows:
\bdesc
  \item[(A)] $m_n=\tlm_n$ for all $n\in\CO{0...N}$.
  \item[(B)] Thus, $\tlM_n = M_n$ for all $n\in\CC{0...N}$.
  \item[(C)] If $m_N=0$ then $\tlm_N=1$. If $m_N=1$ then $\tlm_N=0$.
  \item[(D)] $\tlm_n=0$ for all $n\in\CO{N...N_1}$ (possibly an empty set),
and $\tlm_{N_1}=1$.
  \item[(E)] $m_n=\tlm_n$ for all $n>N_1$.
\edesc
Thus,
\begin{eqnarray}\nonumber
d^N_m &  = &  a_{m+2^N} - a_m 
\quad\eeequals{(*)} \quad
 a^\oo_{k+m+2^N} - a^\oo_{k+m} 
\quad=\quad a^\oo_{\tlM} + a^\oo_{M} \quad \pmod{2}
\\ \nonumber & \eeequals{(\dagger)} &
 \sum_{n=0}^\oo  \lb(\tlm_n\cdot r^n_{\tlM_n} + m_n\cdot r^n_{M_n}\rb)
\quad \eeequals{(ab)}\quad
\sum_{n=N}^\oo 
\lb(\tlm_n\cdot r^n_{\tlM_n} + m_n\cdot r^n_{M_n} \rb)
\\ &\qquad =&
\underbrace{r^N_{M_N}}_{\mathrm{(bc)}}
 \ + \ \
\sum_{n=N+1}^{N_1-1}  m_n \underbrace{r^n_{M_n}}_{{\mathrm{(d)}}}
\ \ + \ \ 
\underbrace{r^{N_1}_{\tlM_{N_1}}}_{{\mathrm{(d)}}}
\ \ + \ \ 
\sum_{n=N_1+1}^\oo  
\underbrace{m_n}_{{\mathrm{(e)}}}
\cdot \lb( r^n_{\tlM_n} + r^n_{M_n}\rb)
\qquad \label{toplitz.irdi.e0} 
\end{eqnarray}
Here, $(*)$ is by equation (\ref{toplitz.irdi.e4}); \ \ 
$(\dagger)$ is by equation (\ref{toplitz.irdi.e5}); \ \ 
(ab) is by (A) and (B); \ \ 
(bc) is by (B) and (C);  \ 
(d) is by (D), \ and
(e) is by (E).

Now, to see {\bf(a)}, let 
\[
S(m) \quad := \quad 
\set{(n,m)}{\maketall \mbox{$r^n_{m}$ appears with
nonzero coefficient in expression 
(\ref{toplitz.irdi.e0})}}.
\]
  In particular, $r^N_{M_N}$ appears in (\ref{toplitz.irdi.e0}), so  $(n_0,m_0):=(N,M_N)$;  thus, $n_0=N$.

To see {\bf(b)}, suppose $m< m'$; hence $m'=m+i$ for some 
$i\in\CO{1...2^N}$.  

Let  $M':=M+i$ and $\tlM':=\tlM+i$.
Suppose $M'=\D\sum_{n=0}^\oo m'_n 2^n$ and 
$\D \tlM'=\sum_{n=0}^\oo \tlm'_n 2^n$.
Define $M'_n$, $\tlM'_n$, and $N_1'$ analogously.
Then, an argument identical to equation (\ref{toplitz.irdi.e0}) yields:
\beqn
\label{toplitz.irdi.e2} 
d^N_{m'}
\quad=\quad
{r^N_{M'_N}}
\ \ + \ \
\sum_{n=N+1}^{N'_1-1}  m'_n {r^n_{M'_n}}
\ \ + \ \ 
{r^{N'_1}_{\tlM'_{N'_1}}}
\ \ + \ \ 
\sum_{n=N'_1+1}^\oo  
{m'_n}
\cdot \lb( r^n_{\tlM'_n} + r^n_{M'_n}\rb)
\eeqn
Now, for all  $n\in\CO{N...\oo}$, \quad 
 $M'_n = M_n+i$ and $\tlM'_n = \tlM_n+i$  (because $i<2^N$);
thus,  $r^n_{M'_n}  =  r^n_{M_n+i}  \not\in\{r^n_{M_n}, r^n_{\tlM_n}\}$ and
$r^n_{\tlM'_n}  =  r^n_{\tlM_n+i}  \not\in\{r^n_{M_n}, r^n_{\tlM_n}\}$.
Thus, every summand of equation (\ref{toplitz.irdi.e2}) is distinct from 
every summand of 
equation (\ref{toplitz.irdi.e0}),  so $S(m')\intsct S(m)=\emptyset$.
\eclaimprf

To see that the random variables
$d^N_{0},\ldots,d^N_{2^N-1}$ are jointly independent, use Claim 1(a):
\[
  d^N_{0}   \  =  \
\sum_{(n,m)\in S(0)} r^{n}_{m},\qquad
  d^N_{1}   \  =  \
\sum_{(n,m)\in S(1)} r^{n}_{m},\quad
 \ldots\ldots\quad
 d^N_{2^N-1}   \  = \  
\sum_{(n,m)\in S(2^N-1)} r^{n}_{m}
\]
The random variables
$\set{r^n_m}{n\in\Natur, \ m \in\CC{1...2^N}}$ are independent,
and Claim 1(b) says 
$S(0),S(1)\ldots,S(2^N-1)$ are pairwise disjoint; thus
$d^N_{0},\ldots,d^N_{2^N-1}$ are jointly independent.

{\bf Lower Decay Rate:} $|\alp|<1$, so if
$N$ is sufficiently large (e.g. $N>L:=  -1/\log_2(\alp)$), then
$\alp^N<\frac{1}{2}$.
Suppose $d^N_m  \  = \   r^{n_0}_{m_0} +  r^{n_1}_{m_1} + \ldots +  r^{n_J}_{m_J}$, as in Claim \ref{toplitz.irdi.c1}(a).  For
all $j\in\CC{0...J}$, let 
$P_j:={\sf Prob}\statement{$\D\sum_{i=j}^J r^{n_i}_{m_i} \ = \ 1$}$.
Thus,
\beq
 \del^N_m\{1\} 
&=&
P_0
\quad\eeequals{(\dagger)}\quad
\rho_{N}\{0\}\cdot P_1 \ + \ \rho_{N}\{1\}\cdot (1-P_1) 
\quad=\quad 
(1-\alp^{N}) \cdot P_1 \ + \ \alp^{N} \cdot (1-P_1)\\
&=&
\alp^{N} \ + \ (1-2\alp^{N})\cdot P_1
\quad\geeeq{(*)}\quad
\alp^{N}
\eeq
$(\dagger)$ is because Claim \ref{toplitz.irdi.c1}(a) says $n_0=N$.\quad
$(*)$ is because $1-2\alp^{N}> 0$, because $\alp^{N} \ < \ \frac{1}{2}$.

{\bf Upper Decay Rate:}  Let
 $K:=\frac{1}{1-\alp}$.  We claim that, for any $N$ and $m$,\quad
  $\del^N_m\{1\} \ \leq \ K\alp^{N}$.

  As before, let $P_j:={\sf Prob}\statement{$\sum_{i=j}^J r^{n_i}_{m_i} \ = \ 1$}$.  For any $j\in\CO{1...J}$, we have
\beqn
P_j
\ \ = \ \ 
(1-\alp^{n_j}) \cdot P_{j+1} \ + \ \alp^{n_j} \cdot (1-P_{j+1})\\
\ \ =\ \  
P_{j+1} \ + \ (1-2P_{j+1})\alp^{n_j}
\ \ \leq\ \ 
P_{j+1} \ + \ \alp^{n_j},\quad
\label{toplitz.irdi.e1}
\eeqn
and $P_J=\alp^{n_J}$.  Hence, 
\[
 \del^N_m\{1\} \quad=\quad  P_0
 \ \ \leeeq{(*)} \ \ 
 \alp^{n_0} + \alp^{n_1} + \ldots + \alp^{n_J}
 \ \ \leq \ \
  \sum_{i=n_0}^\oo \alp^i
 \ \ = \ \
  \frac{\alp^{n_0}}{1-\alp}
 \ \  = \ \
 K\alp^{n_0}
\ \ \eeequals{(\dagger)} \ \ 
 K\alp^{N}.
\]
Here, $(*)$ is obtained by applying equation (\ref{toplitz.irdi.e1}) inductively,
and $(\dagger)$ is because $n_0=N$.
\ethmprf

Thus, if $\frac{1}{\sqrt{2}} < \alp < 1 $, then $\mu$ satisfies the
conditions of Propositions \ref{IRDI.asympt.rand.lemma} 
and \ref{IRDI.entropy.lemma}, so $\mu$ is a zero-entropy, Lucas mixing measure.
Hence, $1+\shift{}$  asymptotically randomizes $\mu$.

{\footnotesize
\bibliographystyle{alpha}
\bibliography{bibliography}
}

\end{document}